\documentclass[12pt]{article}

\usepackage{times}

\usepackage[a4paper,text={128mm,185mm}, centering]{geometry}

\usepackage{changepage}

\usepackage{titlesec}

\titleformat{\section}[hang]
{\bfseries\large}{\thesection.}{1ex}{}

\titleformat{\subsection}[hang]
{\bfseries}{\thesubsection}{1ex}{}

\usepackage{amsthm}
\usepackage{fancyhdr}
\theoremstyle{plain}
\newtheorem{theorem}{Theorem}[section]
\newtheorem{lemma}[theorem]{Lemma}
\newtheorem{proposition}[theorem]{Proposition}
\newtheorem{corollary}[theorem]{Corollary}
\newtheorem{definition}[theorem]{Definition} 

\theoremstyle{definition}
\newtheorem{example}[theorem]{Example}
\newtheorem{remark}[theorem]{Remark}
\newtheorem{notation}[theorem]{Notation}
\newtheorem{terminology}[theorem]{Terminology}

\usepackage{amssymb}
\usepackage{amsmath}
\usepackage{array}
\usepackage{bbm}
\usepackage{bm}
\usepackage{enumerate}
\usepackage[latin1]{inputenc}
\usepackage[T1]{fontenc}
\usepackage{shadow}
\usepackage[all,cmtip,2cell]{xy}
\UseAllTwocells
\usepackage{graphicx}
\usepackage{graphbox}
\usepackage{caption}
\usepackage{subcaption}
\usepackage{float}
\usepackage{multicol}
\usepackage{tikz}
\usepackage{tikz-cd}
\usepackage{longtable}
\usepackage{color}
\usetikzlibrary{matrix,arrows,decorations.pathmorphing}
\usepackage{pdfpages}

\usepackage[colorlinks=true]{hyperref}
\hypersetup{allcolors=[rgb]{0.1,0.1,0.4}}

\newcommand{\CC}{\mathcal{C}}

\newcommand{\sd}{D^s}
\newcommand{\OO}{\mathcal{O}}
\newcommand{\EE}{\mathcal{E}}
\newcommand{\GG}{\mathcal{G}}
\newcommand{\RR}{\mathcal{R}}
\newcommand{\one}{\mathbbm{1}}

\newcommand{\quot}{\delimiter"502F30E\mathopen{}}

\DeclareMathOperator{\Hom}{Hom}
\DeclareMathOperator{\gSet}{gSet}
\DeclareMathOperator{\Set}{Set}
\DeclareMathOperator{\id}{id}

\DeclareMathOperator{\Alg}{Alg}
\DeclareMathOperator{\Comp}{Comp}

\DeclareMathOperator{\Tree}{Tree}
\DeclareMathOperator{\Cell}{Cell}
\DeclareMathOperator{\Psh}{Psh}
\DeclareMathOperator{\cd}{cd}
\DeclareMathOperator{\colim}{colim}
\DeclareMathOperator{\Sesq}{Sesq}

\title{\vskip 5pt  \bf COMPUTADS AND STRING DIAGRAMS FOR $N$-SESQUICATEGORIES}
\author{\itshape\bfseries {Manuel Ara\'{u}jo}}
\date{December 15, 2023}

\begin{document}

\maketitle

\cfoot{}
\thispagestyle{empty}
\vskip 25pt
\begin{adjustwidth}{0.5cm}{0.5cm}
{\small
{\bf R\'esum\'e.} Une $n$-sesquicat\'egorie est un ensemble $n$-globulaire avec des op\'erations de composition strictement associatives et unitaires, qui ne sont cependant pas tenues de satisfaire les lois d'\'echange de Godement qui s'appliquent aux $n$-cat\'egories. Dans \cite{sesquicat}, nous avons montr\'e comment celles-ci peuvent \^etre d\'efinies comme des alg\`ebres sur une monade $T_n^{\sd}$ dont les op\'erations sont des diagrammes de cordes simples. Dans le pr\'{e}sent article, nous donnons une description explicite des polygraphes pour cette monade et nous prouvons que la cat\'egorie associ\'ee de computades est une cat\'egorie de pr\'efaisceaux. Nous utilisons ceci pour d\'ecrire une notation de diagrammes de cordes pour repr\'esenter des compos\'es arbitraires dans des $n$-sesquicat\'egories. Ceci est un pas vers une th\'eorie des diagrammes de cordes pour les $n$-cat\'egories semistrictes.}  \\
{\bf Abstract.} An $n$-sesquicategory is an $n$-globular set with strictly associative and unital composition and whiskering operations, which are however not required to satisfy the Godement interchange laws which hold in $n$-categories. In \cite{sesquicat} we showed how these can be defined as algebras over a monad $T_n^{\sd}$ whose operations are simple string diagrams. In the present paper, we give an explicit description of computads for the monad $T_n^{\sd}$ and we prove that the category of computads for this monad is a presheaf category. We use this to describe a string diagram notation for representing arbitrary composites in $n$-sesquicategories. This is a step towards a theory of string diagrams for semistrict $n$-categories.

\end{adjustwidth}

\section{Introduction}

The use of string diagram notation as a tool for representing composites in higher categories is becoming ever more widespread. This paper is part of a project which aims to give a definition of semistrict $n$-category based on a purely algebraic/combinatorial notion of string diagram. In \cite{sesquicat} we defined a monad $T_n^{\sd}$ on the category of $n$-globular sets, whose operations we call simple string diagrams. We give a generators and relations description of $T_n^{\sd}$, which allows us to characterize its algebras, which we call $n$-sesquicategories, as $n$-globular sets equipped with strictly associative and unital composition and whiskering operations, which however do not satisfy the Godement interchange laws that hold in a strict $n$-category. We think of simple string diagrams as analogous to the globular pasting diagrams used in the definition of the monad $T_n^{str}$ whose algebras are strict $n$-categories (\cite{operads_cats}). In the present paper we study computads for the monad $T_n^{\sd}$ and show how morphisms in an $n$-sesquicategory generated by a computad $C$ can be depicted as general $C$-labelled string diagrams. We also prove that the category of computads for this monad is equivalent to the category of presheaves on a small category of computadic cell shapes. In future work, we will show how to add coherent weak interchange laws to get a notion of semistrict $n$-category,

\subsection{Results}

We now describe the main result in this paper. Denote by $\Comp_{n+1}^{n}$ the category of $(n+1)$-computads for $T_n^{\sd}$, by $\one$ the terminal $(n+1)$-computad and by $F_n(C)$ the free $n$-sesquicategory generated by an $n$-computad $C$. Cells $c\in\one_k$ for $k\leq n+1$ are called  \textbf{$k$-cell shapes} and morphisms $d\in F_n(\one)_k$ for $k\leq n$ are called \textbf{unlabelled $k$-diagrams}. A morphism $x\in F_n(C)$ is said to have shape $d$ if its image in $F_n(\one)$ is $d$. Given such $d$, we construct a computad $\hat{d}$ with the property that $d$-shaped morphisms in a computad $D$ are in canonical bijection with maps $\hat{d}\to D$. Using this, we define a small category $\Cell_{n+1}$ whose objects are cell shapes, together with a fully faithful embedding $\widehat{(-)}:\Cell_{n+1}\hookrightarrow\Comp_{n+1}^{n}$. From this we construct the nerve/realization adjunction \[\begin{tikzcd}
	{|-|:\Psh(\Cell_{n+1}^{n})} & \Comp_{n+1}^{n}:N
	\arrow[""{name=0, anchor=center, inner sep=0}, shift left=2, from=1-1, to=1-2]
	\arrow[""{name=1, anchor=center, inner sep=0}, shift left=2, from=1-2, to=1-1]
	\arrow["\dashv"{anchor=center, rotate=-90}, draw=none, from=0, to=1]
\end{tikzcd}.\]

\begin{theorem}
 
The adjunction $\begin{tikzcd}
	{|-|:\Psh(\Cell_{n+1})} & \Comp_{n+1}^{n}:N
	\arrow[""{name=0, anchor=center, inner sep=0}, shift left=2, from=1-1, to=1-2]
	\arrow[""{name=1, anchor=center, inner sep=0}, shift left=2, from=1-2, to=1-1]
	\arrow["\dashv"{anchor=center, rotate=-90}, draw=none, from=0, to=1]
\end{tikzcd}$ is an equivalence of categories.  
 
\end{theorem}

We now give an outline of the proof. In \cite{sesquicat} we showed that $T_n^{\sd}$ has a presentation with generators $\OO_n$ and relations $\EE_n$. We can describe morphisms in the free $n$-sesquicategory $F_n(C)$ generated by an $n$-computad $C$ as equivalence classes of trees whose internal vertices are labelled by generators in $\OO_n$ and whose leaves are labelled by cells in $C$. The equivalence relation is generated by the relations in $\EE_n$. We then prove that each of these trees has a unique \textbf{normal form} in its equivalence class. This allows us to show that for an unlabelled diagram $d$ the category $\Comp_{n+1}^{n}(d)$ of pairs $(C,x)$, where $C$ is an $(n+1)$-computad and $x$ is a morphism of shape $d$, has an initial object, which we denote $(\hat{d},\tilde{d})$. This allows us to construct the nerve/relization adjunction as mentioned above and then the proof of the Theorem follows by formal arguments from the fact $(\hat{c},\tilde{c})$ is initial, for $c\in\Cell_{n+1}$.

\begin{remark}

In fact our proof of the Theorem above applies to any globular operad presented by generators and relations, as long as this presentation admits a theory of normal forms. See Remark \ref{normal_form_theory} for details.

\end{remark}

The theory of normal forms also provides an \textbf{algorithm} to decide whether two morphisms in the free $n$-sesquicategory $F_n(C)$ generated by $C$ given as composites of generating cells are actually equal.

After we've established this Theorem, we go on to give a description of the diagrammatic interpretation of morphisms in the $n$-sesquicategory generated by $C$ as \textbf{$C$-labelled string diagrams}. Normal forms are an essential ingredient in this description.

\subsection{Related work}

The string diagrammatic calculus for monoidal categories and bicategories is by now well established. Generalizations to Gray $3$-categories also exist, in the theory of surface diagrams (\cite{gray_cats_with_duals}, \cite{hummon_thesis}). Recently there has been a lot of progress in extending this to higher dimensions, with the discovery of the theory of associative $n$-categories (\cite{dorn_thesis}), later developed into the manifold diagrams of \cite{manifold_diagrams}. These manifold diagrams have a combinatorial counterpart, which the authors of \cite{manifold_diagrams} call trusses, which are in turn equivalent to the notion of zigzgags introduced in \cite{high_level_methods} and which forms the basis for an online proof assistant for diagrammatic calculus in higher categories (\cite{homotopyio}). 

There are two main differences between the approach above and the one followed in this paper. The first is that the input of our theory is the simple \textbf{combinatorial} notion of simple string diagram introduced in \cite{sesquicat}, whereas manifold diagrams start from the geometry and obtain from that a combinatorial description, by passing to exit path posets. The second is that we want to produce an \textbf{algebraic} notion of semistrict $n$-categories, by which we mean that these will be algebras over a certain monad on $n$-globular sets. One advantage of the manifold diagrams approach to semistrict $n$-categories is that all coherences are already encoded in the basic cell shapes, whereas we naturally produce a theory of $n$-sesquicategories, to which we then have to add coherent weak interchange laws. The main advantage of our approach is its simplicity, as in a sense everything follows from the combinatorial notion of simple string diagrams introduced in \cite{sesquicat}.

Most closely related to our work is \cite{data_struct_quasistrict}. There the authors develop a framework which is the basis for another proof assistant for  diagrammatic calculus in higher categories (\cite{globular}). The authors have a notion of \textbf{signature}, which corresponds exactly to a computad for $T_n^{\sd}$, and a notion of \textbf{diagram} over a signature, which corresponds exactly to a morphism in the $n$-sesquicatery generated by a computad. In this sense, our work can also be seen as providing a mathematical foundation for the kinds of higher categorical structures implemented by this proof assistant. 

Our work is also related to questions in the general theory of computads (\cite{street_1976}, \cite{street_1987}, \cite{power_1991}, \cite{burroni_1993}, \cite{batanin_computads}, \cite{metayer_2008}). If one considers the monad $T_n^{str}$ whose algebras are strict $n$-categories, then computads consist of presentations for strict $n$-categories. Cells of dimension $k\leq n$ are generating $k$-morphisms and $(n+1)$-cells are relations. The cells of the terminal $n$-computad for $T_n^{str}$ are the most general $n$-categorical cell shapes and the morphisms in the $n$-category generated by it can be thought of as general unlabelled pasting diagrams. One would then like to say that the category of computads for $T_n^{str}$ is a category of presheaves on the cell shapes, but this turns out to be false (\cite{3computads_closed}, \cite{direct_comp_closed}), essentially because of the Eckman-Hilton argument. This lead to the question of finding conditions on monads or restrictions on allowable cells in the associated computads that guarantee that one obtains a presheaf category (\cite{computads_slices_operads},\cite{nonunital_polygraphs},\cite{computads_omega}). Our paper can also be seen as a continuation of this line of research, providing a monad on $n$-globular sets which is related to $T_n^{str}$ and whose category of computads is a presheaf category. On pasting diagrams, see also \cite{forest_2022}.

Finally, the motivation for developing this theory was to be able to use string diagrams to prove results about higher categories. In \cite{fibrations} we develop a string diagram calculus for strict $4$-categories and we use it to prove a result about fibrations of mapping $4$-groupoids. In \cite{adj3} and \cite{adj4} we use this string diagram calculus to prove coherence results for adjunctions in $3$ and $4$-categories. In \cite{araujo_thesis}, we use a string diagram caculus for strict monoidal $3$-categories to prove a coherence result for $3$-dualizable objects in strict symmetric monoidal $3$-categories. 

After the appearance of the present paper on the arXiv, an independent proof of the fact that computads for $n$-sesquicategories form a presheaf category has appered in \cite{forest_mimram}. The authors define $n$-sesquicategories directly by generating operations and relations, so their theory does not mention the combinatorics of simple string diagrams. Their use of rewriting theory to establish the existence of normal forms is a very interesting alternative to our methods in Section 4 of the present paper. To go from normal forms to the main result, they then appeal to Makkai's criterion for presheaf categories. Our approach to this in Section 5 gives a shorter and more direct proof.

\subsection{Future work}

 One can construct a monad $T_n^{ss}$ by adding $(k+1)$-operations (resp. relations) to $T_n^{\sd}$ connecting pairs of simple $k$-string diagrams that map to the same $k$-pasting diagram under the map of monads $T_n^{\sd}\to T_n^{str}$, for $k\leq n-1$ (resp. $k=n$). By constrution, this comes with a contractible map $T_n^{ss}\to T_n^{str}$. We expect the methods in this paper should apply to show that the associated category of $n$-computads is a presheaf category (although the category of $(n+1)$-computads is not, because of necessary interchange relations between $n$-dimensional operations). We can define \textbf{semistrict $n$-categories} as $T_n^{ss}$-algebras. By construction, they will admit a string diagram calculus. Moreover, a result conjectured in \cite[6.2.3]{griffiths_thesis} suggests a possible way of proving that any weak $n$-category is equivalent to a semistrict $n$-category in this sense. This is the subject of ongoing research and we will explore it in future papers.

We are also interested in finding finite descriptions of $T_n^{ss}$. In an upcoming paper, we show how to construct $T_3^{ss}$ by adding a finite set of generators and relations to the monad $T_3^{\sd}$. We then show that its algebras agree with Gray $3$-categories. We are working on extending this to dimension $4$. 

Once the definitions of semistrict $3$ and $4$-categories are in place, we can extend the coherence results for adjunctions of \cite{adj3} and \cite{adj4} to this setting. We will then put this together to extend the coherence result for $3$-dualizable objects of \cite{araujo_thesis} to this setting. An extension of this result to the fully weak setting would give a finite presentation of the framed fully extended $3$-dimensional bordism category, by the Cobordism Hypothesis (\cite{baez_dolan},\cite{lurie},\cite{ayala_francis},\cite{grady_pavlov}).

\section{Background}

Denote by $\gSet_n$ the category of $n$-globular sets. Given a finitary monad $T:\gSet_n\to\gSet_n$ one can define categories $\Comp_k^{T}$ of computads for $T$, for $k=0,\cdots, n+1$, together with adjunctions 

\[\begin{tikzcd}
	{F_k:\Comp_k^T} & \Alg_T:V_k.
	\arrow[""{name=0, anchor=center, inner sep=0}, shift left=2, from=1-1, to=1-2]
	\arrow[""{name=1, anchor=center, inner sep=0}, shift left=2, from=1-2, to=1-1]
	\arrow["\dashv"{anchor=center, rotate=-90}, draw=none, from=0, to=1]
\end{tikzcd}\] This is done inductively, by defining a $k$-computad $C$ to be a tuple $(C_k,C_{\leq k-1},s,t)$ where $C_k$ is a set, which we call the set of $k$-cells of $C$, $C_{\leq k-1}$ is a $(k-1)$-computad, and $s,t:C_k\to F_{k-1}(C_{\leq k-1})_{k-1}$ satisfy the globularity relations $ss=st$ and $ts=tt$. One then defines $F_k$, for $k\leq n$, by the pushout \[\begin{tikzcd}
	{T_n^{\sd}(C_{k}\times\partial\theta^{(k)})} & {T_n^{\sd}(C_{k}\times\theta^{(k)})} \\
	{F_{k-1}(C_{\leq k-1})} & {F_{k}(C),}
	\arrow[from=1-1, to=2-1]
	\arrow[from=1-1, to=1-2]
	\arrow[from=1-2, to=2-2]
	\arrow[from=2-1, to=2-2]
	\arrow["\ulcorner"{anchor=center, pos=0.050, rotate=180}, shift right=2, draw=none, from=2-2, to=1-1]
\end{tikzcd}\] where $\theta^{(k)}$ is the globular set represented by $k$. For $k=n+1$, we replace the inclusion $\partial\theta^{(k)}\hookrightarrow \theta^{(k)}$ by the collpase $\partial\theta^{(n+1)}\to\theta^{(n)}$. Similarly, one defines $V_k$ by a pullback. See \cite{csp_thesis} for a detailed exposition of this theory of computads (for the original references, see \cite{street_1976}, \cite{street_1987}, \cite{power_1991}, \cite{burroni_1993} and \cite{batanin_computads}).

\begin{remark}
 
There are incusion maps $\Comp_k^T\hookrightarrow\Comp_{k+1}^T$ for $k\leq n$, so we can think of $k$-computads as $(n+1)$-computads. For this reason, we sometimes write $\Comp^T$ instead of $\Comp_{n+1}^T$ and use the term computad to refer to an $(n+1)$-computad. 
 
\end{remark}

In \cite{sesquicat} we introduced a monad $T_n^{\sd}$ on globular sets, based on a notion of \textbf{simple string diagram} and we defined an \textbf{$n$-sesquicategory} as an algebra over this monad. The is a map $T_n^{\sd}\to T_n^{str}$ to the monad for strict $n$-categories, so any strict $n$-category is an $n$-sesquicategory. In fact $n$-sesquicategories are just strict $n$-categories without the interchange laws. 

\begin{notation}
 
We denote by $\Sesq_n$ the category of $T_n^{\sd}$-algebras. 
 
\end{notation}

In \cite{sesquicat} we gave a presentation of $T_n^{\sd}$ by generators $\OO_n$ and relations $\EE_n$. There we think of the generators as simple string diagrams, but here we interact with the monad $T_n^{\sd}$ only through this presentation, so we may as well view the generators as symbols. There is a generator $\circ_{i,j}$ for each $i,j=1,\cdots, n$ and a generator $u_i$ for each $i=1,\cdots, n$. Given an $n$-sesquicategory $\CC$, the generator $\circ_{i,j}$ induces a map $\circ_{i,j}^{\CC}:\CC_i\times_{\CC_m}\CC_j\to\CC_M$, where $m=\min\{i,j\}$ and $M=\max\{i,j\}$. We call this \textbf{composition} when $i=j$ and \textbf{whiskering} when $i\neq j$. The generator $u_i$ induces a map $u_i^{\CC}:\CC_{i-1}\to \CC_i$ and we call $u_i^{\CC}(x)$ the \textbf{identity} on $x$. The relations in $\EE_n$ essentially express the associativity and unitality of $\circ_{i,j}$ (these relations also appear below Definition \ref{relation}). 

\begin{notation}
We denote by $\OO_n$ and $\EE_n$ the sets of generators and relations for $T_n^{\sd}$ introduced in \cite{sesquicat} and described in the preceding paragraph.
\end{notation}

\begin{remark}
 
The monad $T_n^{\sd}$ corresponds to an $n$-globular operad $\sd_n$ and the presentation by generators and relations corresponds to a presentation of the globular operad in the sense of \cite{griffiths_2023}.
 
\end{remark}

We also characterize $n$-sesquicategories inductively as categories $\CC$ equipped with a lift of the $\Hom$ functor

\[\xymatrix{ & \Sesq_{n-1}\ar@{}[d]|-{=}\ar[rd]^{(-)_0} &  \\\CC^{op}\times\CC\ar@{.>}[ru]^{\underline{\Hom}_\CC}\ar[rr]_{\Hom_\CC} & & \Set ,}\] but this will not be relevant in the present paper.

We now briefly review the generators and relations description of $T_n^{\sd}$, which our description of computads in the present paper will build on. This discussion will be informal, see \cite{sesquicat} for details.

\begin{definition}
	
Let $X$ be an $n$-graded set. A $k$-dimensional $(\OO_n,X)$-labelled tree is a rooted tree $T$, together with

\begin{enumerate}
	
\item a labelling of its internal vertices $I(T)$ by generators in $\OO_n$;
\item a labelling of its leaves $L(T)$ by elements in $X$;
\item a bijection between the incoming edges at an internal vertex and the inputs of the associated generator;

\end{enumerate}

such that

\begin{enumerate}
	
 \item the root label has dimension $k$;
 \item the source of each incoming edge at an internal vertex has a label of the appropriate dimension.
 
\end{enumerate} 

We denote the set of $k$-dimensional $(\OO_n,X)$-labeled trees by $\Tree_n^{\OO}(X)(k)$ or $\Tree_n^{\OO}(X)_k$.
	
\end{definition}

\begin{terminology}
An \textbf{$(\OO_n,X)$-labelled subtree} (or simply \textbf{subtree} for short) of an $(\OO_n,X)$-labelled tree $T$ consists of all vertices (internal and leaves) that can be reached from a chosen internal vertex of $T$ (the root of the subtree) by travelling towards the leaves.

An \textbf{$\OO_n$-labelled tree} is a tree with a labelling of its vertices by $\OO_n$. An \textbf{$\OO_n$-labelled subtree} of an $(\OO_n,X)$-labelled tree $T$ is a subtree of $T$ in the usual sense, containing  no leaves and inheriting the $\OO_n$-labelling.

\end{terminology}

When $X$ is an $n$-globular set, we can define source and target maps $s,t:\Tree_n^{\OO}(X)_k\to\Tree_n^{\OO}(X)_{k-1}$, although in general they won't satisfy the globularity relation. 

\begin{definition}
 
An \textbf{$n$-preglobular set} is an $n$-graded set $X=\coprod_{i=0}^nX_i$ equipped with source and target maps $s,t:X_k\to X_{k-1}$. A \textbf{globular relation} on $X$ is a relation $\sim$ such that \begin{enumerate}
                                                                                                                                                                                                                              \item if $x\sim\tilde{x}$ then $s(x)\sim s(\tilde{x})$ and $t(x)\sim t(\tilde{x})$;
                                                                                                                                                                                                                                                                                                                                                                                                                                                            \item $ss(x)\sim st(x)$ and $ts(x)\sim tt(x)$.
                                                                                                                                                                                                                                                                                                                                                                                                                                                                                          \end{enumerate} Note that this means the quotient $X\quot\sim$ is an $n$-globular set.

\end{definition}

So given an $n$-globular set $X$, we have an $n$-preglobular set $\Tree_n^{\OO}(X)$.

\begin{definition}
 
We define an $n$-preglobular subset $\Tree_n^{\OO,\EE}(X)\subset\Tree_n^{\OO}(X)$ of $\stackrel{\epsilon}{=}$-compatible trees, equipped with a preglobular relation $\stackrel{\epsilon}{=}$. The definition is by induction on height. The relation $\stackrel{\epsilon}{=}$ is generated by the relations in $\EE_n$. A tree is $\stackrel{\epsilon}{=}$-compatible if for every subtree of the form $x\to\circ_{i,j}\leftarrow y$ we have $s^{i-m+1}(x)\stackrel{\epsilon}{=}t^{j-m+1}(y)$, where $m=\min\{i,j\}$.

\end{definition}

Finally we define $\widetilde{\Tree}_n^{\OO,\EE}(X):=\Tree_n^{\OO,\EE}(X)\quot\stackrel{\epsilon}{=}$ and we show that this defines a monad on $n$-globular sets. We construct a map of monads $$\varphi:\widetilde{\Tree}_n^{\OO,\EE}\to T_n^{\sd}.$$ Each generator in $\OO_n$ corresponds to a simple string diagram, so one can use composition of simple string diagrams to produce this map.

\begin{theorem}[\cite{sesquicat}]
 
The map $\varphi:\widetilde{\Tree}_n^{\OO,\EE}\to T_n^{\sd}$ is an isomorphism of monads.

\end{theorem}

\section{Computads for $T_n^{\sd}$}

We give an explicit description of computads for $T_n^{\sd}$ and of the $n$-sesquicategories generated by them, which we will later show is equivalent to the notion described in the previous section. We will simply call them computads, leaving the monad $T_n^{\sd}$ implicit.

\begin{definition}\label{def_computad}
 
Given $k\leq n+1$, an \textbf{$(n,k)$-precomputad} (or simply \textbf{$k$-precomputad}, leaving $n$ implicit) $C$ consists of sets $C_i$ for $0\leq i\leq k$, together with maps $s,t:C_i\to\Tree_n^{\OO}(C_{\leq i-1})_{i-1}$ for $1\leq i\leq k$. 
 
\end{definition}

In the definition below we use the following notation for grafting of trees.

\begin{notation}

 Given and $(\OO_n,C)$-labelled tree $x\in \Tree_n^{\OO}(C)_{i-1}$ we denote by $$x \to u_i$$ the $(\OO_n,C)$-labelled tree obtained by adding a new new vertex to $x$ labelled by $u_i$ and an edge from the root of $x$ to this new vertex. The new vertex now becomes the root of this new tree. Similarly, given $x\in\Tree_n^{\OO}(C)_i$ and $y\in\Tree_n^{\OO}(C)_j$ we denote by $$x\to \circ_{i,j} \leftarrow y$$ the $(\OO_n,C)$-labelled tree obtained by adding a new root labelled by $\circ_{i,j}$.

\end{notation}

\begin{definition}
 
Given a $k$-precomputad $C$, we define source and target maps $s,t:\Tree_n^{\OO}(C)_{i}\to\Tree_n^{\OO}(C)_{i-1}$, for $1\leq i\leq n$. For trees of height zero, these are the maps $s,t:C_i\to\Tree_n^{\OO}(C)_{i-1}$. For trees of nonzero height, we use the following inductive formulas for $s$, where $j<i$ and $x$ and $y$ are trees with appropriate dimensions in each case. The map $t$ is defined by the same formulas, replacing every instance of $s$ with $t$.

\begin{center}\begin{tabular}{l}

$s(x\to u_i)=x;$ 

\\ \\

$s(x\to \circ_{i,i} \leftarrow y)=s(y);$ 

\\ \\

 $s(x\to \circ_{j,i}\leftarrow y)=x\to \circ_{j,i-1} \leftarrow s(y);$ 
 
 \\ \\

 $s(x\to\circ_{i,j} \leftarrow y)=s(x)\to \circ_{i-1,j} \leftarrow y.$ 
 
 \end{tabular}\end{center}
 
\end{definition}

\begin{remark}\label{height}
 
Since the source or target of a $k$-cell may be an arbitrary $(\OO_n,C)$-labelled tree, the source and target maps above can increase the height of trees. This is in contrast to the situation of \cite{sesquicat}, where we considered $\Tree_n^{\OO}(X)$ for a globular set $X$. However, these maps always decrease the dimension of the tree, so the arguments in \cite{sesquicat} which relied on induction on the height of the tree can now be replaced by simultaneous induction on both the height and the dimension of the tree, as we will do below in Definition \ref{relation}.
 
\end{remark}

The following definitions refer to each other and should be interpreted by mutual induction.

\begin{definition}

Given $k\leq n+1$, an \textbf{$(n,k)$-computad} (or simply \textbf{$k$-computad}, leaving $n$ implicit) $C$ consists of sets $C_i$ for $0\leq i\leq k$, together with maps $s,t:C_i\to\Tree_n^{\OO,\EE}(C_{\leq i-1})_{i-1}$ for $1\leq i\leq k$, such that $ss(x)\stackrel{\epsilon}{=} st(x)$ and $ts(x)\stackrel{\epsilon}{=} tt(x)$ for all $x\in C_i$.

\end{definition}

\begin{terminology}

A \textbf{computad} is an \textbf{$(n,k)$}-computad, where $n$ is usually implicit in the context and $k\leq n+1$ is arbitrary.
 
\end{terminology}

\begin{notation}

Let $X$ be an $n$-graded set. We denote by $$\tau_{\leq h}\Tree_n^{\OO}(X)\subset \Tree_n^{\OO}(X)$$ the $n$-graded subset consisting of trees of height at most $h$. 
 
\end{notation}

The definition that follows is almost identical to the analogous one in \cite{sesquicat}. The only difference is the one explained in Remark~\ref{height}.

\begin{definition}\label{relation} Let $C$ be a computad. For each $k$, we define, by induction on $h$, subsets  $\tau_{\leq h}\Tree_n^{\OO,\EE}(C)_k\subset\tau_{\leq h}\Tree_n^{\OO}(C)_k$ equipped with a relation $\stackrel{\epsilon}{=}_{h}$. Elements in $\tau_{\leq h}\Tree_n^{\OO,\EE}(C)_k$ are called \textbf{$\stackrel{\epsilon}{=}_{h-1}$-compatible}. We say that $x\in\Tree_n^{\OO}(C)_k$ is $\stackrel{\epsilon}{=}$-compatible if it is $\stackrel{\epsilon}{=}_h$-compatible for some $h$ and define $\Tree_n^{\OO,\EE}(C)_k\subset\Tree_n^{\OO}(C)_k$ the set of $\stackrel{\epsilon}{=}$-compatible elements. Finally, we define the relation $\stackrel{\epsilon}{=}$ on $\Tree_n^{\OO,\EE}(C)_k$ by declaring $x\stackrel{\epsilon}{=}\tilde{x}$ when $x\stackrel{\epsilon}{=}_h\tilde{x}$ for some $h$. The definition is by overall induction on $k$ and is presented below.

\end{definition} When $h=0$, we let $\tau_{\leq 0}\Tree_n^{\OO,\EE}(C)_k:=\tau_{\leq 0}\Tree_n^{\OO}(C)_k=C_k$ and the relation $\stackrel{\epsilon}{=}_0$ is $=$. 

Now consider $h\geq 1$. Any $x\in\tau_{\leq h}\Tree_n^{\OO}(C)_k$ of height zero is $\stackrel{\epsilon}{=}_{h-1}$-compatible. Let $x\in\tau_{\leq h-1}\Tree_n^{\OO}(C)_i$, $y\in\tau_{\leq h-1}\Tree_n^{\OO}(C)_j$ and $m=\min\{i,j\}$. Then $$\xymatrixcolsep{1pc}\xymatrixrowsep{1pc}\xymatrix{ & \circ_{i,j} &  \\ x\ar[ru] & & y\ar[lu]}$$  is $\stackrel{\epsilon}{=}_{h-1}$-compatible if and only if $x,y$ are $\stackrel{\epsilon}{=}_{h-2}$-compatible and $s^{i-m+1}(x)\stackrel{\epsilon}{=}t^{j-m+1}(y)$. Moreover, $x\to u_{i+1}$ is $\stackrel{\epsilon}{=}_{h-1}$-compatible if and only if $x$ is $\stackrel{\epsilon}{=}_{h-2}$-compatible. Now we must define the globular relation $\stackrel{\epsilon}{=}_{h}$ on $\tau_{\leq h}\Tree_n^{\OO,\EE}(C)_k$.

If $x,y\in\tau_{\leq h}\Tree_n^{\OO,\EE}(C)_k$ have height zero and $x\stackrel{\epsilon}{=}_{0}y$, then $x\stackrel{\epsilon}{=}_{h}y$. 

Let $i\leq k$, $x\in \tau_{\leq h-2}\Tree_n^{\OO,\EE}(C)_{i-1}$, $y\in \tau_{\leq h-1}\Tree_n^{\OO,\EE}(C)_k$. If $x\stackrel{\epsilon}{=} t^{k-i+1}(y)$, then \begin{center}\begin{tabular}{lccc} $(\lambda_{i,k}):$ & $\vcenter{\vbox{\xymatrixcolsep{1pc}\xymatrixrowsep{1pc}\xymatrix{ & \circ_{i,k} &  \\ u_i\ar[ru] & & y\ar[lu]\\ x\ar[u] & & }}}$ & $\stackrel{\epsilon}{=}_{h}$ & $y$.\end{tabular}\end{center}

Let $i\leq k$, $x\in \tau_{\leq h-1}\Tree_n^{\OO,\EE}(C)_k$, $y\in \tau_{\leq h-2}\Tree_n^{\OO,\EE}(C)_{i-1}$. If $s^{k-i+1}(x)\stackrel{\epsilon}{=}y$, then \begin{center}\begin{tabular}{lccc} $(\rho_{k,i}):$ & $\vcenter{\vbox{\xymatrixcolsep{1pc}\xymatrixrowsep{1pc}\xymatrix{ & \circ_{k,i} &  \\ x\ar[ru] & & u_i\ar[lu]\\  & & y\ar[u]}}}$ & $\stackrel{\epsilon}{=}_{h}$ & $x$.\end{tabular}\end{center}

Let $i<k$, $x\in \tau_{\leq h-2}\Tree_n^{\OO,\EE}(C)_i$, $y\in \tau_{\leq h-2}\Tree_n^{\OO,\EE}(C)_{k-1}$. If $s(x)\stackrel{\epsilon}{=}t^{k-i}(y)$, then \begin{center}\begin{tabular}{lccc} $(\rho_{i,k}):$ & $\vcenter{\vbox{\xymatrixcolsep{1pc}\xymatrixrowsep{1pc}\xymatrix{ & \circ_{i,k} &  \\ x\ar[ru] & & u_k\ar[lu]\\ & & y\ar[u]}}}$ & $\stackrel{\epsilon}{=}_{h}$ & $\vcenter{\vbox{\xymatrixcolsep{1pc}\xymatrixrowsep{1pc}\xymatrix{ & u_k &  \\  & \circ_{i,k-1}\ar[u] & \\ x\ar[ru] & & y\ar[lu]}}}$.\end{tabular}\end{center}

Let $i<k$, $x\in \tau_{\leq h-2}\Tree_n^{\OO,\EE}(C)_{k-1}$, $y\in \tau_{\leq h-2}\Tree_n^{\OO,\EE}(C)_i$. If $s^{k-i}(x)\stackrel{\epsilon}{=}t(y)$, then \begin{center}\begin{tabular}{lccc} $(\lambda_{k,i}):$ & $\vcenter{\vbox{\xymatrixcolsep{1pc}\xymatrixrowsep{1pc}\xymatrix{ & \circ_{k,i} &  \\ u_k\ar[ru] & & y\ar[lu]\\ x\ar[u] & & }}}$ & $\stackrel{\epsilon}{=}_{h}$ & $\vcenter{\vbox{\xymatrixcolsep{1pc}\xymatrixrowsep{1pc}\xymatrix{ & u_k &  \\  & \circ_{k-1,i}\ar[u] & \\ x\ar[ru] & & y\ar[lu]}}}$.\end{tabular}\end{center}

Let $k\geq 1$ and $x,y,z\in \tau_{\leq h-2}\Tree_n^{\OO,\EE}(C)_k$. If $s(x)\stackrel{\epsilon}{=}t(y)$ and $s(y)\stackrel{\epsilon}{=}t(z)$, then \begin{center}\begin{tabular}{lccc} $(\circ_{k,k,k}):$ & $\vcenter{\vbox{\xymatrixcolsep{1pc}\xymatrixrowsep{1pc}\xymatrix{ & & \circ_{k,k} &  \\ & \circ_{k,k}\ar[ru] & & z\ar[lu] \\ x\ar[ru] & & y\ar[lu] & }}}$ & $\stackrel{\epsilon}{=}_{h}$ & $\vcenter{\vbox{\xymatrixcolsep{1pc}\xymatrixrowsep{1pc}\xymatrix{ & \circ_{k,k} & &  \\ x\ar[ru] & & \circ_{k,k}\ar[lu] & \\ & y\ar[ru] & & z\ar[lu]}}}$.\end{tabular}\end{center}

Let $i<k$, $x,y\in \tau_{\leq h-2}\Tree_n^{\OO,\EE}(C)_i$, $z\in \tau_{\leq h-2}\Tree_n^{\OO,\EE}(C)_k$. If $s(x)\stackrel{\epsilon}{=}t(y)$ and $s(y)\stackrel{\epsilon}{=}t^{k-i+1}(z)$, then \begin{center}\begin{tabular}{lccc} $(\circ_{i,i,k}):$ & $\vcenter{\vbox{\xymatrixcolsep{1pc}\xymatrixrowsep{1pc}\xymatrix{ & & \circ_{i,k} &  \\ & \circ_{i,i}\ar[ru] & & z\ar[lu] \\ x\ar[ru] & & y\ar[lu] & }}}$ & $\stackrel{\epsilon}{=}_{h}$ & $\vcenter{\vbox{\xymatrixcolsep{1pc}\xymatrixrowsep{1pc}\xymatrix{ & \circ_{i,k} & &  \\ x\ar[ru] & & \circ_{i,k}\ar[lu] & \\ & y\ar[ru] & & z\ar[lu]}}}$.\end{tabular}\end{center}

Let $i<k$, $x,z\in \tau_{\leq h-2}\Tree_n^{\OO,\EE}(C)_i$, $y\in \tau_{\leq h-2}\Tree_n^{\OO,\EE}(C)_k$. If $s(x)\stackrel{\epsilon}{=}t^{k-i+1}(y)$ and $s^{k-i+1}(y)\stackrel{\epsilon}{=}t(z)$, then \begin{center}\begin{tabular}{lccc} $(\circ_{i,k,i}):$ & $\vcenter{\vbox{\xymatrixcolsep{1pc}\xymatrixrowsep{1pc}\xymatrix{ & & \circ_{k,i} &  \\ & \circ_{i,k}\ar[ru] & & z\ar[lu] \\ x\ar[ru] & & y\ar[lu] & }}}$ & $\stackrel{\epsilon}{=}_{h}$ & $\vcenter{\vbox{\xymatrixcolsep{1pc}\xymatrixrowsep{1pc}\xymatrix{ & \circ_{i,k} & &  \\ x\ar[ru] & & \circ_{k,i}\ar[lu] & \\ & y\ar[ru] & & z\ar[lu]}}}$.\end{tabular}\end{center}

Let $i<k$, $x\in \tau_{\leq h-2}\Tree_n^{\OO,\EE}(C)_k$, $y,z\in \tau_{\leq h-2}\Tree_n^{\OO,\EE}(C)_i$. If $s^{k-i+1}(x)\stackrel{\epsilon}{=}t(y)$ and $s(y)\stackrel{\epsilon}{=}t(z)$, then \begin{center}\begin{tabular}{lccc} $(\circ_{k,i,i}):$ & $\vcenter{\vbox{\xymatrixcolsep{1pc}\xymatrixrowsep{1pc}\xymatrix{ & & \circ_{k,i} &  \\ & \circ_{k,i}\ar[ru] & & z\ar[lu] \\ x\ar[ru] & & y\ar[lu] & }}}$ & $\stackrel{\epsilon}{=}_{h}$ & $\vcenter{\vbox{\xymatrixcolsep{1pc}\xymatrixrowsep{1pc}\xymatrix{ & \circ_{k,i} & &  \\ x\ar[ru] & & \circ_{i,i}\ar[lu] & \\ & y\ar[ru] & & z\ar[lu]}}}$.\end{tabular}\end{center}

Let $i<k$, $x\in \tau_{\leq h-2}\Tree_n^{\OO,\EE}(C)_i$, $y,z\in \tau_{\leq h-2}\Tree_n^{\OO,\EE}(C)_k$. If $s(x)\stackrel{\epsilon}{=}t^{k-i+1}(y)$ and $s(y)\stackrel{\epsilon}{=}t(z)$, then \begin{center}\begin{tabular}{lccc} $(\circ_{i,k,k}):$ & $\vcenter{\vbox{\xymatrixcolsep{1pc}\xymatrixrowsep{1pc}\xymatrix{ & \circ_{i,k} & &  \\ x\ar[ru] & & \circ_{k,k}\ar[lu] & \\ & y\ar[ru] & & z\ar[lu]}}}$ & $\stackrel{\epsilon}{=}_{h}$ & $\vcenter{\vbox{\xymatrixcolsep{1pc}\xymatrixrowsep{1pc}\xymatrix{ & & \circ_{k,k} & & \\ & \circ_{i,k}\ar[ru] & & \circ_{i,k}\ar[lu] & \\ x\ar[ru] & & y\ar[lu]\text{  }x\ar[ru] & & z\ar[lu]}}}$.\end{tabular}\end{center}

Let $i<k$, $x,y\in \tau_{\leq h-2}\Tree_n^{\OO,\EE}(C)_k$, $z\in \tau_{\leq h-2}\Tree_n^{\OO,\EE}(C)_i$. If $s(x)\stackrel{\epsilon}{=}t(y)$ and $s^{k-i+1}(y)\stackrel{\epsilon}{=}t(z)$, then \begin{center}\begin{tabular}{lccc} $(\circ_{k,k,i}):$ & $\vcenter{\vbox{\xymatrixcolsep{1pc}\xymatrixrowsep{1pc}\xymatrix{ & & \circ_{k,i} &  \\ & \circ_{k,k}\ar[ru] & & z\ar[lu] \\ x\ar[ru] & & y\ar[lu] & }}}$ & $\stackrel{\epsilon}{=}_{h}$ & $\vcenter{\vbox{\xymatrixcolsep{1pc}\xymatrixrowsep{1pc}\xymatrix{ & & \circ_{k,k} & & \\ & \circ_{k,i}\ar[ru] & & \circ_{k,i}\ar[lu] & \\ x\ar[ru] & & z\ar[lu]\text{  }y\ar[ru] & & z\ar[lu]}}}$.\end{tabular}\end{center}

Let $i<j<k$ and take $x\in \tau_{\leq h-2}\Tree_n^{\OO,\EE}(C)_i$, $y\in \tau_{\leq h-2}\Tree_n^{\OO,\EE}(C)_j$ and $z\in \tau_{\leq h-2}\Tree_n^{\OO,\EE}(C)_k$. If $s(x)\stackrel{\epsilon}{=}t^{j-i+1}(y)$ and $s(y)\stackrel{\epsilon}{=}t^{k-j+1}(z)$, then \begin{center}\begin{tabular}{lccc} $(\circ_{i,j,k}):$ & $\vcenter{\vbox{\xymatrixcolsep{1pc}\xymatrixrowsep{1pc}\xymatrix{ & \circ_{i,k} & &  \\ x\ar[ru] & & \circ_{j,k}\ar[lu] & \\ & y\ar[ru] & & z\ar[lu]}}}$ & $\stackrel{\epsilon}{=}_{h}$ & $\vcenter{\vbox{\xymatrixcolsep{1pc}\xymatrixrowsep{1pc}\xymatrix{ & & \circ_{j,k} & & \\ & \circ_{i,j}\ar[ru] & & \circ_{i,k}\ar[lu] & \\ x\ar[ru] & & y\ar[lu]\text{  }x\ar[ru] & & z\ar[lu]}}}$.\end{tabular}\end{center}

Let $i<j<k$ and take $x\in \tau_{\leq h-2}\Tree_n^{\OO,\EE}(C)_i$, $y\in \tau_{\leq h-2}\Tree_n^{\OO,\EE}(C)_k$ and $z\in \tau_{\leq h-2}\Tree_n^{\OO,\EE}(C)_j$. If $s(x)\stackrel{\epsilon}{=}t^{k-i+1}(y)$ and $s^{k-j+1}(y)\stackrel{\epsilon}{=}t(z)$, then \begin{center}\begin{tabular}{lccc} $(\circ_{i,k,j}):$ & $\vcenter{\vbox{\xymatrixcolsep{1pc}\xymatrixrowsep{1pc}\xymatrix{ & \circ_{i,k} & &  \\ x\ar[ru] & & \circ_{k,j}\ar[lu] & \\ & y\ar[ru] & & z\ar[lu]}}}$ & $\stackrel{\epsilon}{=}_{h}$ & $\vcenter{\vbox{\xymatrixcolsep{1pc}\xymatrixrowsep{1pc}\xymatrix{ & & \circ_{k,j} & & \\ & \circ_{i,k}\ar[ru] & & \circ_{i,j}\ar[lu] & \\ x\ar[ru] & & y\ar[lu]\text{  }x\ar[ru] & & z\ar[lu]}}}$.\end{tabular}\end{center}

Let $i<j<k$ and take $x\in \tau_{\leq h-2}\Tree_n^{\OO,\EE}(C)_j$, $y\in \tau_{\leq h-2}\Tree_n^{\OO,\EE}(C)_k$ and $z\in \tau_{\leq h-2}\Tree_n^{\OO,\EE}(C)_i$. If $s(x)\stackrel{\epsilon}{=}t^{k-j+1}(y)$ and $s^{k-i+1}(y)\stackrel{\epsilon}{=}t(z)$, then \begin{center}\begin{tabular}{lccc} $(\circ_{j,k,i}):$ & $\vcenter{\vbox{\xymatrixcolsep{1pc}\xymatrixrowsep{1pc}\xymatrix{ & & \circ_{k,i} &  \\ & \circ_{j,k}\ar[ru] & & z\ar[lu] \\ x\ar[ru] & & y\ar[lu] & }}}$ & $\stackrel{\epsilon}{=}_{h}$ & $\vcenter{\vbox{\xymatrixcolsep{1pc}\xymatrixrowsep{1pc}\xymatrix{ & & \circ_{j,k} & & \\ & \circ_{j,i}\ar[ru] & & \circ_{k,i}\ar[lu] & \\ x\ar[ru] & & z\ar[lu]\text{  }y\ar[ru] & & z\ar[lu]}}}$.\end{tabular}\end{center}

Let $i<j<k$ and take $x\in \tau_{\leq h-2}\Tree_n^{\OO,\EE}(C)_k$, $y\in \tau_{\leq h-2}\Tree_n^{\OO,\EE}(C)_j$ and $z\in \tau_{\leq h-2}\Tree_n^{\OO,\EE}(C)_i$. If $s^{k-j+1}(x)\stackrel{\epsilon}{=}t(y)$ and $s^{j-i+1}(y)\stackrel{\epsilon}{=}t(z)$, then \begin{center}\begin{tabular}{lccc} $(\circ_{k,j,i}):$ & $\vcenter{\vbox{\xymatrixcolsep{1pc}\xymatrixrowsep{1pc}\xymatrix{ & & \circ_{k,i} &  \\ & \circ_{k,j}\ar[ru] & & z\ar[lu] \\ x\ar[ru] & & y\ar[lu] & }}}$ & $\stackrel{\epsilon}{=}_{h}$ & $\vcenter{\vbox{\xymatrixcolsep{1pc}\xymatrixrowsep{1pc}\xymatrix{ & & \circ_{k,j} & & \\ & \circ_{k,i}\ar[ru] & & \circ_{j,i}\ar[lu] & \\ x\ar[ru] & & z\ar[lu]\text{  }y\ar[ru] & & z\ar[lu]}}}$.\end{tabular}\end{center}

Let $x,\tilde{x}\in \tau_{\leq h-1}\Tree_n^{\OO,\EE}(C)_{k-1}$. If $x\stackrel{\epsilon}{=}_{h-1} \tilde{x}$, then \begin{center}\begin{tabular}{lccc} $(u_k):$ & $\vcenter{\vbox{\xymatrixcolsep{1pc}\xymatrixrowsep{1pc}\xymatrix{u_k \\ x\ar[u]}}}$ & $\stackrel{\epsilon}{=}_{h}$ &  $\vcenter{\vbox{\xymatrixcolsep{1pc}\xymatrixrowsep{1pc}\xymatrix{u_k \\ \tilde{x}\ar[u]}}}$.\end{tabular}\end{center}

Let $x,\tilde{x}\in \tau_{\leq h-1}\Tree_n^{\OO,\EE}(C)_i$, $y,\tilde{y}\in\tau_{\leq h-1}\Tree_n^{\OO,\EE}(C)_j$ and $m=\min\{i,j\}$. If $x\stackrel{\epsilon}{=}_{h-1} \tilde{x}$, $y\stackrel{\epsilon}{=}_{h-1} \tilde{y}$, $s^{i-m+1}(x)\stackrel{\epsilon}{=} t^{j-m+1}(y)$ and $s^{i-m+1}(\tilde{x})\stackrel{\epsilon}{=} t^{j-m+1}(\tilde{y})$ then \begin{center}\begin{tabular}{lccc} $(\circ_{i,j}):$ & $\vcenter{\vbox{\xymatrixcolsep{1pc}\xymatrixrowsep{1pc}\xymatrix{ & \circ_{i,j} &  \\ x\ar[ru] & & y\ar[lu]}}}$ & $\stackrel{\epsilon}{=}_{h}$ & $\vcenter{\vbox{\xymatrixcolsep{1pc}\xymatrixrowsep{1pc}\xymatrix{ & \circ_{i,j} &  \\ \tilde{x}\ar[ru] & & \tilde{y}\ar[lu]}}}$.\end{tabular}\end{center}  

\begin{lemma}
 
The construction above defines an $n$-preglobular subset $\Tree_n^{\OO,\EE}(C)\subset\Tree_n^{\OO}(C)$ with a globular relation $\stackrel{\epsilon}{=}$, meaning we have

\begin{enumerate}
 \item if $x\in\Tree_n^{\OO}(C)$ is $\stackrel{\epsilon}{=}$-compatible, then so are $s(x)$ and $t(x)$;
 
 \item if $x\stackrel{\epsilon}{=}\tilde{x}$ then $s(x)\stackrel{\epsilon}{=}s(\tilde{x})$ and $t(x)\stackrel{\epsilon}{=}t(\tilde{x})$;
 
 \item if $x$ is $\stackrel{\epsilon}{=}$-compatible, then $ss(x)\stackrel{\epsilon}{=}st(x)$ and $ts(x)\stackrel{\epsilon}{=}tt(x)$.
\end{enumerate}

\end{lemma}

\begin{proof}
 
The proof is very similar to the one for the analogous result in \cite{sesquicat}, the only difference being the one already mentioned in Remark~\ref{height}. \end{proof}

\begin{definition}

Let $0\leq k\leq n$ and let $C$ be a $k$-computad. We write $$\overline{\Tree}_n^{\OO,\EE}(C):=\Tree_n^{\OO,\EE}(C)\quot \stackrel{\epsilon}{=}.$$
 
\end{definition}

\begin{definition}\label{relc}

Let $C$ be an $(n+1)$-computad. We define a relation $\stackrel{\epsilon,C}{=}$ on $\Tree_n^{\OO,\EE}(C)$ by adding the new equation $s(x)\stackrel{\epsilon,C}{=}t(x)$ for each $x\in C_{n+1}$.
 
\end{definition}

\begin{lemma}
 
Let $C$ be an $(n+1)$-computad. Then $\stackrel{\epsilon,C}{=}$ is a globular relation on $\Tree_n^{\OO,\EE}(C)$.
 
\end{lemma}

\begin{proof}
 
This is easy to check. \end{proof}
 
\begin{definition}

Let $C$ be an $(n+1)$-computad. We write $$\overline{\Tree}_n^{\OO,\EE}(C):=\Tree_n^{\OO,\EE}(C)\quot \stackrel{\epsilon,C}{=}.$$
 
\end{definition}

\begin{remark}
 
The fact that $\stackrel{\epsilon}{=}$ and $\stackrel{\epsilon,C}{=}$ are globular relations implies $\overline{\Tree}_n^{\OO,\EE}(C)$ is an $n$-globular set. Using the isomorphism of monads $\widetilde{\Tree}_n^{\OO,\EE}\to T^{\sd}_n$ allows us to define a $T_n^{\sd}$ action on $\overline{\Tree}_n^{\OO,\EE}(C)$ by simply grafting trees. We refer to $\overline{\Tree}_n^{\OO,\EE}(C)$ as the \textbf{$n$-sesquicategory presented by $C$}. When $C$ is an $n$-computad this is a \textbf{free $n$-sesquicategory}. When $C$ is an $(n+1)$-computad, this is a quotient of the free $n$-sesquicaegory generated by $C_{\leq n}$ by the relations in $C_{n+1}$.
 
\end{remark}

\begin{definition}
 
Given $k$-computads $C,D$ a map $f:C\to D$ is a collection of maps $f_i: C_i\to D_i$ such that $s(f_i(x))\stackrel{\epsilon}{=}f_{i-1}(s(x))$ and $t(f_i(x))\stackrel{\epsilon}{=}f_{i-1}(t(x))$ for all $x\in C_i$, where we have inductively used the map on trees induced by a map of $(k-1)$-computads.

A map $f:C\to D$ induces a map $f:\Tree_n^{\OO,\EE}(C)\to\Tree_n^{\OO,\EE}(D)$ by applying $f$ to leaf labels.
 
\end{definition}

\begin{definition}
 
For $k\leq n+1$, we denote by $\Comp_k^n$ the category of $(n,k)$-computads and $(n,k)$-computad maps.
 
\end{definition}

\begin{remark}
 
Adding empty sets of cells provides an inclusion map $$\Comp_k^n\hookrightarrow\Comp_{k+1}^n,$$ for $k\leq n$, so we can think of $k$-computads as $(n+1)$-computads. For this reason, we sometimes write $\Comp^n$ instead of $\Comp^n_{n+1}$ and refer to $(n+1)$-computads simply as \textbf{computads}. In fact we will ususally denote this category simply by $\Comp$, leaving $n$ implicit.
 
\end{remark}

\begin{lemma}\label{pushout1}
 
Let $C$ be a computad. Then the following diagram is a pushout, for $k\leq n$.

\[\begin{tikzcd}
	{T_n^{\sd}(C_{k}\times\partial\theta^{(k)})} & {T_n^{\sd}(C_{k}\times\theta^{(k)})} \\
	{\overline{\Tree}_n^{\OO,\EE}(C_{\leq k-1})} & {\overline{\Tree}_n^{\OO,\EE}(C_{\leq k})}
	\arrow[from=1-1, to=1-2]
	\arrow[from=1-1, to=2-1]
	\arrow[from=1-2, to=2-2]
	\arrow[from=2-1, to=2-2]
	\arrow["\ulcorner"{anchor=center, pos=0, rotate=180}, shift right=2, draw=none, from=2-2, to=1-1]
\end{tikzcd}\]
 
\end{lemma}

\begin{proof}
 
We must show that functors $$\varphi_{\leq k}:{\overline{\Tree}_n^{\OO,\EE}(C_{\leq k})}\to\CC$$ correspond to pairs $(\varphi_{\leq k-1},\varphi_{k})$, where $\varphi_{\leq k-1}:{\overline{\Tree}_n^{\OO,\EE}(C_{\leq k-1})}\to\CC$ is a functor and $\varphi_{k}:C_{k}\to\CC_{k}$ is a map, such that $\varphi_{\leq k-1}(s(x)) = s(\varphi_{k}(x))$ and $\varphi_{\leq k-1}(t(x)) = t(\varphi_{k}(x))$ for all $x\in C_{k}$. This is clear.\end{proof}

\begin{lemma}\label{pushout2}

Let $C$ be a computad. Then the following diagram is a pushout.

\[\begin{tikzcd}
	{T_n^{\sd}(C_{n+1}\times\partial\theta^{(n+1)})} & {T_n^{\sd}(C_{n+1}\times\theta^{(n)})} \\
	{\overline{\Tree}_n^{\OO,\EE}(C_{\leq n})} & {\overline{\Tree}_n^{\OO,\EE}(C)}
	\arrow[from=1-1, to=1-2]
	\arrow[from=1-1, to=2-1]
	\arrow[from=1-2, to=2-2]
	\arrow[from=2-1, to=2-2]
	\arrow["\ulcorner"{anchor=center, pos=0, rotate=180}, shift right=2, draw=none, from=2-2, to=1-1]
\end{tikzcd}\]
 
\end{lemma}

\begin{proof}

We must show that functors $$\varphi:{\overline{\Tree}_n^{\OO,\EE}(C)}\to\CC$$ correspond to functors $\varphi_{\leq n}:{\overline{\Tree}_n^{\OO,\EE}(C_{\leq n})}\to\CC$ such that $\varphi_{\leq n}(s(x)) = \varphi_{\leq n}(t(x))$ for all $x\in C_{n+1}$. This is clear. \end{proof}

\begin{proposition}\label{computads}
 
For $k\leq n+1$, the canonical map $\Comp^n_k\to\Comp_k^{T_n^{\sd}}$ is an equivalence of categories, and the following diagram commutes up to canonical natural isomorphism. \[\begin{tikzcd}
	{\Comp^n_k} && {\Comp_k^{T_n^{\sd}}} \\
	& {\Sesq_n}
	\arrow["{\overline{\Tree}^{\OO,\EE}_n}"', from=1-1, to=2-2]
	\arrow["{F_k}", from=1-3, to=2-2]
	\arrow["\simeq", from=1-1, to=1-3]
\end{tikzcd}\]
 
\end{proposition}

\begin{proof}
 
Using induction on $k$ and Lemmas~\ref{pushout1} and \ref{pushout2} we get a canonical map $$\Comp^n_k\to\Comp_k^{T_n^{\sd}}$$ such that the above diagram commutes up to canonical natural isomorphism. 

To construct an inverse, given a computad $C\in \Comp_k^{T_n^{\sd}}$ and using induction on $k$ one can view its source and target maps as $$s,t:C_i\to\overline{\Tree}_n^{\OO,\EE}(C_{\leq i-1})_{i-1}$$ for $i\leq k$. Using the axiom of choice to obtain a section $$\overline{\Tree}_n^{\OO,\EE}(C_{\leq i-1})_{i-1}\to\Tree_n^{\OO,\EE}(C_{\leq i-1})_{i-1}$$ of the quotient map,  we finally obtain maps $s,t:C_i\to\Tree_n^{\OO.\EE}(C_{\leq i-1})_{i-1}$ as in Definition \ref{def_computad}. \end{proof}

\begin{remark}
 
One can avoid using the axiom of choice by using instead normal forms, which give an explicit unique representative of each equivalence class in $\overline{\Tree}_n^{\OO,\EE}(C)$.
 
\end{remark}

\section{Normal form}

In this section, given an $n$-computad $C$, we introduce the notion of normal form for elements of $\Tree_n^{\OO,\EE}(C)$. Denoting by $N(C)\subset \Tree_n^{\OO,\EE}(C)$ the $n$-graded subset of elements in normal form, we prove that for any $x\in\Tree_n^{\OO,\EE}(C)$ there exists a unique $n(x)\in N(C)$ such that $n(x)\stackrel{\epsilon}{=}x$.

\begin{remark}
 
The proof below actually gives an algorithm for finding the normal form $n(x)$ associated to any term $x\in\Tree_n^{\OO,\EE}(C)$. Thus it gives an algorithm for deciding whether two such terms are equivalent.
 
\end{remark}

\begin{remark}

If $C$ is an $(n+1)$-computad, every term $x\in\Tree_n^{\OO,\EE}(C)$ still has a unique normal form $n(x)\stackrel{\epsilon}{=}x$. However, any nontrivial relation in $C_{n+1}$ will provide terms $x,y$ such that $x\stackrel{\epsilon,C}{=}y$ and $n(x)\neq n(y)$ (recall Definition \ref{relc}). So normal forms apply most naturally to $n$-computads.
 
\end{remark}

\begin{notation}
 
We write $m(\circ_{i,j}):=\min\{i,j\}$. When $v$ is an internal vertex in an $(\OO_n,C)$-labelled tree with label $\circ_{i,j}$, we write $m(v)=m(\circ_{i,j})$. 
 
\end{notation}

\begin{definition}
 
An $x\in\Tree_n^{\OO}(C)$ is \textbf{$m$-ordered} if for every edge of the form $v\to w$, where $v,w$ are $\circ$-labelled, we have $m(v)<m(w)$.
 
\end{definition}

\begin{definition}
 
An $\OO_n$-labelled tree is \textbf{$m$-constant} if there are no $u$-labelled vertices and for every edge $v\to w$ we have $m(v)=m(w)$.
 
\end{definition}

\begin{definition}

An \textbf{$m$-constant component} of $x\in\Tree_n^{\OO}(C)$ is a maximal $m$-constant $\OO_n$-labelled subtree. 
 
\end{definition}

\begin{definition}

An $m$-constant $\OO_n$-labelled tree is in \textbf{normal form} if it is of one of the following forms:

\begin{enumerate}

\item $\vcenter{\vbox{\xymatrixcolsep{0.5pc}\xymatrixrowsep{0.1pc}\xymatrix{ \circ_{i,k} & & & \\ & \cdots\ar[lu] & &\\ & & \circ_{i,k}\ar[lu] & \\ & & & \circ_{k,i}\ar[lu] \\ & & \cdots\ar[ru] & \\ & \circ_{k,i}\ar[ru] & & }}}$  with $i<k$;

\item $\vcenter{\vbox{\xymatrixcolsep{0.5pc}\xymatrixrowsep{0.1pc}\xymatrix{ & & \circ_{k,k}\\ & \cdots\ar[ru] & \\ \circ_{k,k}\ar[ru] & & }}}$.
 
\end{enumerate}

\end{definition}

\begin{definition}
 
An $x\in\Tree_n^{\OO}(C)$ is in \textbf{normal form} if it is $m$-ordered, it contains no edges of the form $u\to\circ$ and each of its $m$-constant components is in normal form. 
 
\end{definition}

\begin{definition}
 
Let $x\in\Tree_n^{\OO}(C)$. Define its \textbf{cell dimension} to be the maximum of the dimensions of the cells labelling the leaves of $x$. Denote this by $\cd(x)$. 

\end{definition}

\begin{lemma}
 
Let $x,\tilde{x}\in\Tree_n^{\OO,\EE}(C)$ and suppose $x\stackrel{\epsilon}{=}\tilde{x}$. Then $\cd(x)=\cd(\tilde{x})$. 
\end{lemma}

\begin{proof}This is clear.\end{proof}

\begin{lemma}
 
Given $x\in\Tree_n^{\OO,\EE}(C)$, there exists $n(x)\stackrel{\epsilon}{=} x$ which is in normal form. 
 
\end{lemma}

\begin{proof}
 
One uses the defining equations of $\stackrel{\epsilon}{=}$ to rearrange generators. The $(\lambda)$ and $(\rho)$ relations allow us eliminate all units $u_i$ for $i\leq\cd(x)$ and push the other units towards the root. Then the $(\circ)$ relations allow us to pass to an $m$-ordered tree and finally to put each $m$-constant component in normal form. \end{proof}

Now we need to show that this normal form is unique.

\begin{proposition}\label{unique}
 
Let $x,\tilde{x}\in\Tree_n^{\OO,\EE}(C)$ be in normal form and suppose  $x\stackrel{\epsilon}{=}\tilde{x}$. Then $x=\tilde{x}$. 
 
\end{proposition}

We will prove this below. First we reduce to diagrams without $u$-labelled vertices.

\begin{lemma}\label{circ}

Let $x,\tilde{x}\in\Tree_n^{\OO,\EE}(C)_k$ be in normal form and suppose $x\stackrel{\epsilon}{=}\tilde{x}$. Let $\cd:=\cd(x)=\cd(\tilde{x})$. Then $x=(x_{\circ}\to u_{\cd+1}\to\cdots\to u_k)$ and $\tilde{x}=(\tilde{x}_\circ\to u_{\cd+1}\to\cdots\to u_k)$, where $x_{\circ}$ and $\tilde{x}_\circ$ are in normal form, have no $u$-labelled vertices and $x_{\circ}\stackrel{\epsilon}{=}\tilde{x}_{\circ}$.
 
\end{lemma}

\begin{proof}
 
It is obvious that one can decompose elements in normal form into a unit chain and a component containing no units. The only thing that requires proof is the fact that $x_{\circ}\stackrel{\epsilon}{=}\tilde{x}_{\circ}$. This follows from the observation that  $x_\circ=s^{k-\cd}(x)$ and $\tilde{x}_\circ=s^{k-\cd}(\tilde{x})$. \end{proof}

The above Lemma allows us to reduce the proof of Proposition~\ref{unique} to the case where $x$, $\tilde{x}$ have no units. Now we would like to reduce to the case where one gets from $x$ to $\tilde{x}$ without introducing units along the way.

\begin{definition}
 
Define $\Tree_n^{\OO(\circ)}(C)\subset\Tree_n^{\OO}(C)$ to be the preglobular subset consisting of those trees not containing any $u$-labelled vertices. We then define a preglobular subset $\Tree_n^{\OO(\circ),\EE(\circ)}(C)\subset\Tree_n^{\OO(\circ)}(C)$ of $\stackrel{\circ}{=}$-compatible trees with a globular relation $\stackrel{\circ}{=}$, in exactly the same way we defined $\stackrel{\epsilon}{=}$ and $\stackrel{\epsilon}{=}$-compatibility, except we omit all equations involving $u$.

\end{definition}

\begin{definition}
 
We define the \textbf{reduction} $r(x)\in\Tree_n^{\OO(\circ)}$ of $x\in\Tree_n^{\OO}(C)$ inductively, as follows. We let $r(x)=x$ when $x$ has height zero. Then, for $i<k$, we let $r(x\to u_k)=\emptyset$ and $$r(x\rightarrow\circ_{k,k}\leftarrow y)=\begin{cases}r(y) & r(x)=\emptyset; \\ r(x) & r(y)=\emptyset; \\ r(x)\rightarrow\circ_{k,k}\leftarrow r(y) & \text{otherwise;}\end{cases}$$ 

$$r(x\rightarrow\circ_{i,k}\leftarrow y)=\begin{cases}r(y) & r(x)=\emptyset; \\ \emptyset & r(y)=\emptyset; \\ r(x)\rightarrow\circ_{i,k}\leftarrow r(y) & \text{otherwise;}\end{cases}$$

$$r(x\rightarrow\circ_{k,i}\leftarrow y)=\begin{cases}r(x) & r(y)=\emptyset; \\ \emptyset & r(x)=\emptyset; \\ r(x)\rightarrow\circ_{k,i}\leftarrow r(y) & \text{otherwise.}\end{cases}$$
 
\end{definition}

\begin{lemma}\label{idemp}
 
If $x\in\Tree_n^{\OO(\circ)}$, then $r(x)=x$. 
 
\end{lemma}

\begin{proof}
 
This is obvious. \end{proof}

\begin{lemma}\label{empty}
 
If $w\in\Tree_n^{\OO,\EE}(C)$ and $r(w)=\emptyset$ then $s(w)\stackrel{\epsilon}{=}t(w)$.  
 
\end{lemma}

\begin{proof}
 
The proof is by induction on the height of $w$. There are four cases, corresponding to the four possible root labels: $u_k$, $\circ_{k,k}$, $\circ_{i,k}$ and $\circ_{k,i}$, for $i<k$. Each of these follows by a simple argument. \end{proof}

\begin{lemma}\label{r_props}
 
Given $w,\tilde{w}\in\Tree_n^{\OO,\EE}(C)$, we have \begin{enumerate}
                                                     \item $r(w)\in\Tree_n^{\OO(\circ),\EE(\circ)}(C)$;
                                                     \item if $r(w)\neq\emptyset$, then $sr(w)\stackrel{\circ}{=}rs(w)$ and $tr(w)\stackrel{\circ}{=}rt(w)$;
                                                     \item if $w\stackrel{\epsilon}{=}\tilde{w}$, then $r(w)\stackrel{\circ}{=}r(\tilde{w})$.
\end{enumerate}
 
\end{lemma}

\begin{proof}
 
The proof is by mutual induction on dimension and height. For \textit{1.} there is  a case for each possible root label of $w$: $u_k$, $\circ_{k,k}$, $\circ_{i,k}$ and $\circ_{k,i}$ ($i<k$). Each of these follows from a simple inductive argument. 

For \textit{2.} there are cases for root labels $u_k$, $\circ_{k,k}$, $\circ_{k-1,k}$, $\circ_{k,k-1}$, $\circ_{i,k}$ and $\circ_{k,i}$ ($i<k-1$). We explain the $\circ_{k,k}$ case and leave the others to the reader. Let $w=(x \to {\circ_{k,k}} \leftarrow y).$ Now $$sr(w)=\begin{cases} sr(x) & r(y)=\emptyset \\ sr(y) & r(y)\neq\emptyset 
\end{cases}$$ and $rs(w)=rs(y)$. If $r(y)\neq\emptyset$, we have $rs(y)\stackrel{\circ}{=}sr(y)$ by induction, so $rs(w)\stackrel{\circ}{=}sr(w)$. When $r(y)=\emptyset$, we need to show that $sr(x)\stackrel{\circ}{=}rs(y)$. We have $r(x)\neq\emptyset$, because $r(w)\neq\emptyset$. Then $sr(x)\stackrel{\circ}{=}rs(x)$ by induction. We also have $s(y)\stackrel{\epsilon}{=}t(y)$ by Lemma~\ref{empty}. Since $w$ is $\stackrel{\epsilon}{=}$-compatible, we have $s(x)\stackrel{\epsilon}{=}t(y)$, so we have $s(x)\stackrel{\epsilon}{=}s(y)$ and then using \textit{3.} by induction we have $rs(x)\stackrel{\circ}{=}rs(y)$, so $sr(x)\stackrel{\circ}{=}rs(y)$.
 
To prove \textit{3.}, there is one case for each of the defining equations of $\stackrel{\epsilon}{=}$. We explain the $(\circ_{k,k,k})$ case, leaving the others to the reader. Let $w$ and $\tilde{w}$ be the left and right hand sides of this equation, respectively. If at least one of the trees $r(x)$, $r(y)$, $r(z)$ is empty, then we get $r(w)=r(\tilde{w})$ and we are done. So we may assume they are all nonempty. In this case we get $r(w)\stackrel{\circ}{=}r(\tilde{w})$ by the same $(\circ_{k,k,k})$ equation, as long as $r(x)$, $r(y)$, $r(z)$ are $\stackrel{\circ}{=}$-compatible, $sr(x)\stackrel{\circ}{=}tr(y)$ and $sr(y)\stackrel{\circ}{=}tr(z)$. This first condition follows from \textit{1.} by induction on height. The second condition follows from \textit{2.} and \textit{3.} by induction on height and dimension. \end{proof}

%
%

\begin{notation}
 
Given $x\in\Tree_n^{\OO}(C)$, we denote by $L(x)$ its set of leaves. Given $\ell \in L(x)$, we denote by $|\ell|$ the dimension of the cell labellng $\ell$. We denote by $L_{\geq i}(x)\subset L(x)$ the set of leaves $\ell$ such that $|\ell|\geq i$.
 
\end{notation}

\begin{definition}
 
Given $x\in\Tree_n^{\OO(\circ)}(C)$, we define $M(x)=\max\{j:|L_{\geq j}(x)|\geq 2\}$. If $x$ only has one leaf, then $M(x)=-\infty$.
 
\end{definition}

\begin{lemma}
 
If $x\stackrel{\circ}{=}\tilde{x}$, then $M(x)=M(\tilde{x})$. 
 
\end{lemma}

\begin{proof}
One just needs to check that this holds for each of the defining equations of $\stackrel{\circ}{=}$, which is easy.\end{proof}

\begin{definition}
 
Given $w\in\Tree_n^{\OO}(C)$, we define a linear ordering of $L(w)$ as follows. When $w=(x\to u_k)$ then $L(w)=L(x)$ and we can just use induction on height. When $w=(x \to \circ_{i,j} \leftarrow y)$ then $L(w)=L(x)\coprod L(y)$ and we define the linear order on $L(w)$ by using the linear orders on $L(x)$ and $L(y)$ provided by induction on height, together with the rule that $\ell_x>\ell_y$ for any $\ell_x\in L(x)$ and $\ell_y\in L(y)$.
 
\end{definition}

\begin{lemma}
 
Let $w,\tilde{w}\in\Tree_n^{\OO(\circ),\EE(\circ)}(C)$ and suppose $w\stackrel{\circ}{=}\tilde{w}$. Let $M:=M(w)=M(\tilde{w})$. Then there is a (necessarily unique) order preserving isomorphism $L_{\geq M}(w)\to L_{\geq M}(\tilde{w})$.
 
\end{lemma}

\begin{proof}
 
One just needs to check this for each of the equations defining $\stackrel{\circ}{=}$. The only equations requring some consideration are $(\circ_{i,k,k})$, $(\circ_{k,k,i})$, $(\circ_{i,j,k})$, $(\circ_{i,k,j})$, $(\circ_{j,k,i})$ and $(\circ_{k,j,i})$, which double some of the leaves. In each case, one can see that this doubling does not affect leaves in $L_{\geq M}$. For example, the $(\circ_{i,k,k})$ equation doubles the leaves in the subtree $x$. But since we have no $u$-labelled vertices, the subtrees $y$ and $z$ must both have at least one leaf labelled by a $k$-cell, so that $w$ must have at least two leaves labelled by $k$-cells, so that $M=k$. Then $i<k$ implies $i<M$, so there are no leaves labelled by cells of dimension $\geq M$ in $x$.  \end{proof}

%
%
%
%
%

\begin{definition}
 
Let $w\in\Tree_n^{\OO(\circ),\EE(\circ)}(C)_k$, let $ M(w)\leq M\leq k$, and let $\ell\in L_{\geq M}(w)$. We define $\sigma^M_{\ell}(w)\in\Tree_n^{\OO(\circ)}(C)_{|\ell|}$ by induction on height as follows. If $w$ has height zero, then it consists of a single leaf $\ell$, and we let $$\sigma^M_{\ell}(\ell)=\ell.$$ If $w$ has nonzero height, then we have a case for each possible root label. For $p,q\in\{M,k\}$ (with at least one equal to $k$) and $i<M$, we let  \begin{longtable}{lll} $\sigma^M_{\ell}(x\to \circ_{p,q}\leftarrow y)$ & $=$ & $\begin{cases} \sigma^M_{\ell}(x), & \ell\in L(x); \\ \sigma^M_{\ell}(y), & \ell\in L(y); \end{cases}$

\\ \\

$\sigma^M_{\ell}(x\to \circ_{i,k}\leftarrow y)$ & $=$ & $(x\to \circ_{i,|\ell|}\leftarrow \sigma^M_{\ell}(y));$

\\ \\

$\sigma^M_{\ell}(x\to \circ_{k,i}\leftarrow y)$ & $=$ & $(\sigma^M_{\ell}(x)\to \circ_{|\ell|,i}\leftarrow y).$

\end{longtable}
 
\end{definition}

\begin{lemma}\label{slices}
 
Let $w\in\Tree_n^{\OO(\circ),\EE(\circ)}(C)_k$, let $M(w)\leq M \leq k$ and let $\ell\in L_{\geq M}(w)$. Then

\begin{enumerate}
 \item $\sigma^M_{\ell}(w)$ is $\stackrel{\circ}{=}$-compatible;
 \item $s^{|\ell|-i+1}(\sigma^M_{\ell}(w))\stackrel{\circ}{=}s^{k-i+1}(w)$ and $t^{|\ell|-i+1}(\sigma^M_{\ell}(w))\stackrel{\circ}{=}t^{k-i+1}(w)$ for $i<M$;
 \item if $w\stackrel{\circ}{=}\tilde{w}$ and $\tilde{\ell}\in L_{\geq M}(\tilde{w})$ is the image of $\ell$, then $\sigma^M_{\ell}(w)\stackrel{\circ}{=}\sigma^M_{\tilde{\ell}}(\tilde{w})$.

\end{enumerate}
 
\end{lemma}

\begin{proof}
 
The proof is by mutual induction on the height of $w$. One proves \textit{1.} easily by splitting into the cases wich appear in the definition of $\sigma^M_{\ell}$ and using \textit{2.} on trees of smaller height.

To prove \textit{2.}, we again split into the cases appearing in the definition of $\sigma^M_{\ell}$. We explain only the case $w=(x\to \circ_{p,q}\leftarrow y)$, as the others are simpler. We also do only $s$, as $t$ is completely analogous. So we compute \begin{longtable}{ll}
                                                                                                                                                                                                                                                                  
$s^{|\ell|-i+1}(\sigma_{\ell}^M(x\to \circ_{p,q}\leftarrow y))$ & $=\begin{cases} s^{|\ell|-i+1}(\sigma_{\ell}^M(x)), & \ell\in L(x) \\ s^{|\ell|-i+1}(\sigma_{\ell}^M(y)), & \ell\in L(y)                                                                                                                                                                                                                                                         \end{cases}$

\\ \\

& $\stackrel{\circ}{=} \begin{cases} s^{p-i+1}(x), & \ell\in L(x) \\ s^{q-i+1}(y), & \ell\in L(y),                                                                                                                                                                                                                                                         \end{cases}$                                                                                                                                                                                                                                                                  \end{longtable} where we used induction. On the other hand $$s^{k-i+1}(x\to \circ_{p,q}\leftarrow y)=s^{q-i+1}(y).$$ Now recall that $\stackrel{\circ}{=}$ is a globular relation, so $ss\stackrel{\circ}{=}st$. Moreover, we have $s^{p-m+1}(x)\stackrel{\circ}{=}t^{q-m+1}(y)$, where $m=\min\{p,q\}$, because $w$ is $\stackrel{\circ}{=}$-compatible. This allows us to compute $$s^{p-i+1}(x)=s^{m-i}s^{p-m+1}(x)\stackrel{\circ}{=}s^{m-i}t^{q-m+1}(y)\stackrel{\circ}{=}s^{m-i}s^{q-m+1}(y)\stackrel{\circ}{=}s^{q-i+1}(y)$$ so we are done.

For \textit{3.} there is one case for each of the defining equations of $\stackrel{\circ}{=}$. The arguments are simple in every case, so we leave them to the reader.
 
\end{proof}

\begin{definition}
 
Let $w\in\Tree_n^{\OO(\circ),\EE(\circ)}(C)$. We define $\sigma_{\ell}(w):=\sigma^{M(w)}_{\ell}(w)$.
 
\end{definition}

\begin{definition}
 
Let $x\in\Tree_n^{\OO(\circ)}(C)$. We define $$H(x)=|\{m(v):v\in I(T)\}|.$$
 
\end{definition}

\begin{lemma}
 
Let $x,\tilde{x}\in\Tree_n^{\OO(\circ),\EE(\circ)}(C)$ and suppose $x\stackrel{\circ}{=}\tilde{x}$. Then $$H(x)=H(\tilde{x}).$$
 
\end{lemma}

\begin{proof}

One checks this is true for each of the equations defining $\stackrel{\circ}{=}$, which is easy. \end{proof}

\begin{lemma}\label{nounits}
 
Let $x,\tilde{x}\in\Tree_n^{\OO(\circ),\EE(\circ)}(C)$ be in normal form and suppose $x\stackrel{\circ}{=}\tilde{x}$. Then $x=\tilde{x}$.  
 
\end{lemma}

\begin{proof}
 
The proof is by induction on $H:=H(x)=H(\tilde{x})$. If $H=0$, then $x$, $\tilde{x}$ both have height zero, so they must be equal as the relation $\stackrel{\circ}{=}$ is just $=$ on elements of height zero. 

Now suppose $H\geq 1$ and let $l=|L_{\geq M}(x)|=|L_{\geq M}(\tilde{x})|$. Then $x$, $\tilde{x}$ each consist of a maximal $m$-constant component containing the root, which we denote $x_0$ and $\tilde{x}_0$, to which are grafted trees $x_1,\cdots,x_l$ and $\tilde{x}_1,\cdots,\tilde{x_l}$. Moreover, it is easy to see that $x_i=\sigma_{\ell_i}(x)$ and $\tilde{x}_i=\sigma_{\tilde{\ell}_i}(\tilde{x})$, where $L_{\geq M}(x)=\{\ell_1<\cdots<\ell_l\}$ and $L_{\geq M}(\tilde{x})=\{\tilde{\ell}_1<\cdots<\tilde{\ell}_l\}$.  By Lemma~\ref{slices}, we must have $\sigma_{\ell_i}(x)\stackrel{\circ}{=}\sigma_{\tilde{\ell}_i}(\tilde{x})$ and so by induction we have $\sigma_{\ell_i}(x)=\sigma_{\tilde{\ell}_i}(\tilde{x})$. 

Now we must show $x_0=\tilde{x}_0$. If $M=k$ then both must be equal to $\circ_{k,k}\to\cdots\to\circ_{k,k}$, where there are $l-1$ copies of $\circ_{k,k}$. If $M<k$, then only one of the leaves $\ell_i$ will be labelled by a $k$-cell, let it be $\ell_p$. This also means $\tilde{\ell}_p$ must be the only leaf in $\tilde{x}$ labelled by a $k$-cell. Then both $x_0$ and $\tilde{x}_0$ must be equal to $\circ_{k,M}\to\cdots\to\circ_{k,M}\to\circ_{M,k}\to\cdots\to\circ_{M,k}$ where we have $p-1$ copies of $\circ_{k,M}$ and $l-p$ copies of $\circ_{M,k}$. \end{proof}

\begin{proof}[Proof of Proposition~\ref{unique}]
 
We have $x,\tilde{x}\in\Tree_n^{\OO,\EE}(C)_k$, both in normal form, and $x\stackrel{\epsilon}{=}\tilde{x}$. By Lemma~\ref{circ}, we can assume that $x,\tilde{x}$ contain no $u$-labelled vertices. By Lemmas~\ref{idemp} and \ref{r_props}, we then have $x\stackrel{\circ}{=}\tilde{x}$ and we can apply Lemma~\ref{nounits} to conclude $x=\tilde{x}$.\end{proof}

%

\begin{notation}

We denote by $N(C)\subset\Tree_n^{\OO,\EE}(C)$ the $n$-graded subset consisting of terms in normal form.
 
\end{notation}

\begin{corollary}
 
Let $C$ be an $n$-computad. Then the map $$n:\overline{\Tree}_n^{\OO,\EE}(C)\to N(C)$$ sending an equivalence class to its unique representative in normal from is an inverse to the map $N(C)\hookrightarrow \Tree_n^{\OO,\EE}(C) \to \overline{\Tree}_n^{\OO,\EE}(C)$. 
 
\end{corollary}

\begin{proof}
 
This follows directly from Proposition \ref{unique}.\end{proof}

We record here two important properties of normal forms.

\begin{lemma}\label{normal_form}
 
If $x$ is in normal form then any subtree of $x$ is also in normal form. If $\phi$ is a map of computads, then $\phi(x)$ is in normal form if and only if $x$ is. 
 
\end{lemma}

\begin{proof} This is clear. \end{proof}

%
%

\section{$\Comp$ is a presheaf category}

We now define a category $\Cell_{n+1}$ of computadic cell shapes of dimension $\leq n+1$, with a fully faithful functor $\Cell_{n+1}\hookrightarrow\Comp^n_{n+1}$, which we ususally denote $\Cell\hookrightarrow\Comp$. Then we construct the associated nerve/realization adjunction \[\begin{tikzcd}
	{|-|:\Psh(\Cell)} & \Comp:N
	\arrow[""{name=0, anchor=center, inner sep=0}, shift left=2, from=1-1, to=1-2]
	\arrow[""{name=1, anchor=center, inner sep=0}, shift left=2, from=1-2, to=1-1]
	\arrow["\dashv"{anchor=center, rotate=-90}, draw=none, from=0, to=1]
\end{tikzcd}\] and prove that it is an equivalence of categories. 

The main ingredient is Proposition \ref{initial}, which uses normal forms in an essential way. The rest of the section could probably be shortened by appealing to the theory of familial representability (\cite{carboni_johnstone_1995},\cite{nonunital_polygraphs}). We choose to present the arguments here for the reader's convenience, as this does not take too much space.

\begin{definition}
 
Denote by $\one$ the terminal computad. One can define it inductively by saying that it has exactly one $0$-cell and exactly one $k$-cell $x$ with $s(x)=x_0$ and $t(x)=x_1$ for each ordered pair $(x_0,x_1)$ of parallel $(k-1)$-morphisms in  $F_{k-1}(\one_{\leq k-1})$. 
 
\end{definition}

\begin{notation}
 
For any computad $C$, we denote by $\sigma:C\to\one$ the unique map of computads. 
 
\end{notation}

We think of the cells in $\one$ as \textbf{cell shapes} and of the morphisms in the free $n$-sesquicategory $F_n(\one)$ as \textbf{unlabelled diagrams}. Now we show how one can associate a computad to each unlabelled diagram.

\begin{definition}
 
Let $k\leq n$ and let $d\in F_n(\one)$ be an unlabelled $k$-diagram. We define a category $\Comp(d)$ as follows. Its objects are pairs $(C,x)$, where $C$ is a computad, $x\in F_n(C)$ and $F_n(\sigma)(x)=d$. A morphism $(C,x)\to (D,y)$ is a map of computads $\varphi:C\to D$ such that $F_n(\varphi)(x)=y$.
 
\end{definition}

Now we show that the category $\Comp(d)$ has an initial object. For this we will need to take colimits in $\Comp$. The following result seems to be well known, but not having found a suitable reference we include a simple proof here.

\begin{lemma}\label{formal}
 
Let $T:\gSet_n\to\gSet_n$ be a finitary monad and let $k\leq n+1$. For each $m\leq k$, denote by $[-]_m:\Comp_k^T\to\Set$ the functor taking a $k$-computad to its set of $m$-cells. Then the following hold:

\begin{enumerate}
 \item the category $\Comp_k^T$ is cocomplete;
 \item each functor $[-]_m$ is cocontinuous;
 \item the functors $[-]_m$ for $m=0,\cdots,k$ jointly reflect isomorphisms.
 
\end{enumerate}

\end{lemma}

\begin{proof}
 
Let $\Gamma:I\to\Comp_k^T$ be a diagram. To prove \textit{1.} and \textit{2.} we may as well assume $m=k$, otherwise we can pass to the underlying diagram of $m$-computads. We construct, by induction on $k$, a $k$-computad $C$ which will be the colimit of this diagram. Define its set of $k$-cells to be $C_k:=\colim_i [\Gamma(i)]_k$ and its underlying $(k-1)$-computad as $C_{\leq k-1}:=\colim_i [\Gamma(i)]_{\leq k-1}$. Now define source and target maps $$s,t:\colim_i [\Gamma(i)]_k\to [F_{k-1}(\colim_i [\Gamma(i)]_{\leq k-1})]_{k-1}$$ by the composite \begin{center}\begin{tabular}{l}$\Gamma(i)_k\to [F_{k-1}(\Gamma(i)_{\leq k-1})]_{k-1}\to\colim_i [F_{k-1}(\Gamma(i)_{\leq k-1})]_{k-1}\to$  \\ \\ $\to[\colim_iF_{k-1}(\Gamma(i)_{\leq k-1})]_{k-1}=[F_{k-1}(\colim_i \Gamma(i)_{\leq k-1})]_{k-1}$\end{tabular}\end{center} where the equality comes from the fact that $F_{k-1}$ is left adjoint and the last arrow is induced by the maps $$[F_{k-1}([\Gamma(i)]_{\leq k-1})]_{k-1}\to[\colim_iF_{k-1}([\Gamma(i)]_{\leq k-1})]_{k-1}$$ on sets of $(k-1)$-morphisms associated to the canonical maps of $T$-algebras $F_{k-1}([\Gamma(i)]_{\leq k-1})\to\colim_iF_{k-1}([\Gamma(i)]_{\leq k-1})$. Now one needs to check that $s,t$ satisfy globularity and that the construction has the right universal property. This is straightforward. Point \textit{3.} is easy to prove by induction. \end{proof}

\begin{proposition}\label{initial}
 
For each unlabelled diagram $d\in F_n(\one)$ the category $\Comp(d)$ has an initial object, which we denote $(\hat{d},\tilde{d})$.
 
\end{proposition}

\begin{proof}
 
We construct $(\hat{d},\tilde{d})$ by induction on the dimension of $d$ and on the height of its normal form. 

If $d$ is a $0$-diagram then it consists of a single $0$-cell. Then $\hat{d}$ is the $0$-computad with a single $0$-cell and $\tilde{d}$ is the diagram consisting of that $0$-cell. Now suppose $d$ consists of a single $k$-cell. By induction on dimension and the fact that $ss(d)=st(d)$ and $ts(d)=tt(d)$, we have the following diagram.

\[\begin{tikzcd}
	& {\widehat{s(d)}} \\
	{\widehat{s^2(d)}} && {\widehat{t^2(d)}} \\
	& {\widehat{t(d)}}
	\arrow[hook, from=2-1, to=1-2]
	\arrow[hook', from=2-3, to=1-2]
	\arrow[hook', from=2-1, to=3-2]
	\arrow[hook, from=2-3, to=3-2]
\end{tikzcd}\]

We build $\hat{d}$ by taking the colimit of this diagram in $\Comp$ and then adding a $k$-cell $\tilde{d}:\widetilde{s(d)}\to\widetilde{t(d)}$. It's now easy to see, by induction on dimension, that $(\hat{d},\tilde{d})$ is an initial object in $\Comp(d)$.

Now suppose $d$ has normal form $x \to {\circ_{i,j}} \leftarrow y.$ Let $m=\min\{i,j\}$ and let $x\cap y=s^{i-m+1}(x)=t^{j-m+1}(y)$. By induction on height and dimension, we have a diagram \[\begin{tikzcd}
	{\widehat{x\cap y}} & {\hat{x}} \\
	{\hat{y}}
	\arrow[hook', from=1-1, to=2-1]
	\arrow[hook, from=1-1, to=1-2]
\end{tikzcd}.\] We let $\hat{d}$ be the pushout of this diagram in $\Comp$ and $$\tilde{d}=(\tilde{x} \to \circ_{i,j} \leftarrow \tilde{y}),$$ where we take $\tilde{x},\tilde{y}$ in normal form. Since $d$ is in normal form, so are $x$ and $y$, by Lemma \ref{normal_form}. Therefore, again by Lemma \ref{normal_form} and uniqueness of normal forms, $\tilde{x},\tilde{y}$ map to $x,y$ in $\Tree_n^{\OO,\EE}(\one)$. Then $(\tilde{x} \to \circ_{i,j} \leftarrow \tilde{y})$ maps to $(x \to {\circ_{i,j}} \leftarrow y)$ in $\Tree_n^{\OO,\EE}(\one)$ and therefore it is in normal form, by Lemma \ref{normal_form}. Given $(C,m)\in\Comp(d)$, a map $(\hat{d},\tilde{d})\to(C,m)$ is given by maps $f:\hat{x}\to C$ and $g:\widehat{y}\to C$ such that $$(f(\tilde{x})\to\circ_{i,j}\leftarrow g(\tilde{y}))\stackrel{\epsilon}{=}m.$$ Note that the left hand side is already in normal form, by Lemma \ref{normal_form}. This means the normal form of $m$ must be equal to this, by uniqueness of normal form. This determines $f(\tilde{x}),g(\tilde{y})\in F_n(C)$ uniquely, because it determines their normal forms as the two evident subtrees of the normal form of $m$. Then by induction this determines $f,g$ uniquely, so we conclude that $(\hat{d},\tilde{d})$ is initial. 

Finally suppose $d$ has normal form $x\to u_i$. Then take $\hat{d}=\hat{x}$ and $\tilde{d}=\tilde{x}\to u_i$.\end{proof}

\begin{remark}

The pair $(\hat{d},\tilde{d})$ corresponds to what is called a \textbf{polyplex} in \cite{burroni_2012} and \cite{nonunital_polygraphs}. In \cite{nonunital_polygraphs}, the essential condition for establishing that a certain class of polygraphs forms a presheaf category is the fact that the groups of autmorphisms of polyplexes are trivial.
 
\end{remark}

\begin{remark}\label{no_normal_forms}
 
In fact, we don't need normal forms to construct these computads. We only need them to prove that they are initial. One can construct $(\hat{x},\tilde{x})$  for any term $x\in\Tree_n^{\OO,\EE}(\one_{\leq n})$ by the same inductive procedure used above. When $x\stackrel{\epsilon}{=}y$ is one of the generating equations in $\EE_n$, we obtain an isomorphism $\varphi:\hat{x}\to\hat{y}$ such that $\varphi(\tilde{x})\stackrel{\epsilon}{=}\tilde{y}$ by the same generating equation. This means $(\hat{d},\tilde{d})$ is well defined up to isomorphism. If $x$ is in normal form and $\varphi$ is an automorphism of $\hat{x}$ such that $\varphi(\tilde{x})\stackrel{\epsilon}{=}\tilde{x}$, then $\varphi(\tilde{x})=\tilde{x}$, as they are both in normal form. This implies $\varphi=\id$. So in the presence of normal forms there are no automorphisms, so $(\hat{d},\tilde{d})$ is well defined up to unique isomorphism and it is initial in $\Comp(d)$.
 
\end{remark}

\begin{example}
 
It is well known (\cite{3computads_closed},\cite{direct_comp_closed}) that, for $n\geq 2$, the category of $3$-computads for the monad $T_n^{str}$  whose algebras are strict $n$-categories, is not a presehaf category. This example illustrates why the the above Proposition fails in this case. Denote by $\one$ the terminal computad for $T_n^{str}$ and let $s:\id_*\Rightarrow\id_*$ be the unique $2$-cell in $\one$ whose source and target are the identity on the unique $0$-cell. We construct a diagram $d\in F_2(\one)_2$ consisting of the vertical composite $s\circ s$. Consider $(\hat{d},\tilde{d})\in\Comp^{T_n^{str}}(d)$ defined by letting $\hat{d}$ be the computad consisting of a $0$-cell $*$, together with two $2$-cells $\alpha,\beta:\id_*\Rightarrow\id_*$, and $\tilde{d}$ the vertical composite $\alpha\circ\beta$. By the Eckmann-Hilton argument, we have $\alpha\circ\beta=\beta\circ\alpha$, so $\hat{d}$ admits a nontrivial automorphism which maps $\tilde{d}\mapsto\tilde{d}$, namely the one that permutes $\alpha$ and $\beta$. If $\Comp(d)$ had an initial object $I$, then the unique map $I\to (\hat{d},\tilde{d})$ would be invariant under composition with this automorphism. This would mean that $\alpha,\beta$ are not in the image of the map, so $I$ contains only $0$-cells, which is absurd.

In order to show that $\Comp$ is a presheaf category, what we actually need is the fact that $\Comp(c)$ has an initial object when $c\in\one_k$ is a computadic cell shape. This will fail for any $3$-cell shape whose source or target is the diagram $d$ above. 
 
\end{example}

\begin{definition}
 
Let $c\in\one_{n+1}$ be an $(n+1)$-cell shape. We define $\Comp(c)$ to be the category of pais $(C,x)$ where $C$ is a computad and $x\in C_{n+1}$ is an $(n+1)$-cell such that $\sigma(x)=c$. 
 
\end{definition}

\begin{corollary}\label{cellinitial}

For each $k\leq n+1$ and each $k$-cell $c\in\one_k$ the category $\Comp(c)$ has an initial object, denoted $(\hat{c},\tilde{c})$.
 
\end{corollary}

\begin{proof}
 
For $k\leq n$, this is just Proposition \ref{initial}. For an $(n+1)$-cell $c$, we use Proposition \ref{initial} to construct the diagram \[\begin{tikzcd}
	& {\widehat{s(c)}} \\
	{\widehat{s^2(c)}} && {\widehat{t^2(c)}} \\
	& {\widehat{t(c)}}
	\arrow[hook, from=2-1, to=1-2]
	\arrow[hook', from=2-3, to=1-2]
	\arrow[hook', from=2-1, to=3-2]
	\arrow[hook, from=2-3, to=3-2]
\end{tikzcd}\] and then we take the colimit and add an $(n+1)$-cell $\tilde{c}:\widetilde{s(c)}\to\widetilde{t(c)}$.
\end{proof}

\begin{remark}
 
The pair $(\hat{c},\tilde{c})$ corresponds to what is called a \textbf{plex} in \cite{burroni_2012} and \cite{nonunital_polygraphs} or a \textbf{computope} in \cite{makkai_word}. 
 
\end{remark}

\begin{definition}
 
Let $\Cell_{n+1}$ be the category whose objects are cell shapes $c\in\one_k$ for $k\leq n+1$, and where a morphism $c\to d$ is a map of computads $\hat{c}\to\hat{d}$. We usually denote this simply by $\Cell$. It comes with a fully faithful functor $$\widehat{(-)}:\Cell\hookrightarrow\Comp.$$
 
\end{definition}

\begin{definition}
 
We define the nerve functor $$N:\Comp\to\Psh(\Cell)$$ as the composite $\Comp\hookrightarrow \Psh(\Comp) \to \Psh(\Cell)$ of the Yoneda embedding with the restriction along $\widehat{(-)}$.

\end{definition}

\begin{definition}
 
We define the realization functor by the following left Kan extension, which exists because $\Comp$ is cocomplete. \[\begin{tikzcd}
	\Cell & \Comp \\
	{\Psh(\Cell)}
	\arrow["{\widehat{(-)}}", hook, from=1-1, to=1-2]
	\arrow["{\mathcal{Y}}"', hook', from=1-1, to=2-1]
	\arrow["{|-|}"', dashed, from=2-1, to=1-2]
\end{tikzcd}\] 
 
\end{definition}

 We thus obtain the usual nerve/realization adjunction

    \[\begin{tikzcd}
	{|-|:\Psh(\Cell)} & \Comp:N
	\arrow[""{name=0, anchor=center, inner sep=0}, shift left=2, from=1-1, to=1-2]
	\arrow[""{name=1, anchor=center, inner sep=0}, shift left=2, from=1-2, to=1-1]
	\arrow["\dashv"{anchor=center, rotate=-90}, draw=none, from=0, to=1]
\end{tikzcd}.\]

\begin{theorem}
 
The adjunction $$\begin{tikzcd}
	{|-|:\Psh(\Cell_{n+1})} & \Comp_{n+1}^n:N
	\arrow[""{name=0, anchor=center, inner sep=0}, shift left=2, from=1-1, to=1-2]
	\arrow[""{name=1, anchor=center, inner sep=0}, shift left=2, from=1-2, to=1-1]
	\arrow["\dashv"{anchor=center, rotate=-90}, draw=none, from=0, to=1]
\end{tikzcd}$$ is an equivalence.  
 
\end{theorem}

\begin{proof}
 
By \cite{computads_omega}[Proposition 5.14] it is enough to show that the functors $$\Comp(\widehat{c},-):\Comp\to\Set,$$ for $c\in\Cell$, are cocontinuous and jointly reflect isomorphisms. This follows easily from Lemma~\ref{formal} and the fact that for $c\in\one_k$ we have $$\Comp(\widehat{c},C)=\{x\in C_k : \sigma(x)=c\},$$ which follows from Corrollary \ref{cellinitial}. \end{proof}

\begin{remark}\label{normal_form_theory}
 
All results in this section hold, with the same proofs, for any $n$-globular operad given by generators and relations as long as it admits a suitable theory of normal forms. More precisely, given a presentation $(\GG,\RR)$ for an $n$-globular operad, what we need is an $n$-graded subset $N(C)\subset \Tree_n^{\GG,\RR}(C)$ of terms in normal form, for each $n$-computad $C$, with the following properties:
\begin{enumerate}
 \item the induced map $N(C)\to\overline{\Tree}_n^{\GG,\RR}(C)$ is an $n$-graded bijection (i.e. there is a unique term in normal form in each equivalence class);
 
 \item each subtree of a tree in normal from is in normal form;
 
 \item given a map of $n$-computads $\phi:C\to D$, we have $\phi(x)\in N(D)$ if and only if $x\in N(C)$. 

\end{enumerate}
 
\end{remark}

\begin{remark}

Because of condition $3.$ in the previous Remark, it is enough to define $N(\one_{\leq n})$ for the terminal computad $\one_{\leq n}$ and then let $x\in N(C)$ if and only if $\sigma(x)\in N(\one_{\leq n})$. It is also enough to check condition $2.$ for the terminal computad. However, it is not enough to check $1.$ for $\one_{\leq n}$, as one can have $x\neq y$ in $N(C)$ such that $\sigma(x)=\sigma(y)$.
 
\end{remark}

\begin{remark}

We now describe an alternative approach to showing that the category of computads for an $n$-globular operad presented by generators and relations is a presheaf category. This was inspired by a discussion with Samuel Mimram about the theory of rewriting.
 
Denote by $\Gamma_n^{\GG,\RR}$ the free groupoid on the graph with vertices the terms $x\in\Tree_n^{\GG,\RR}(\one_{\leq n})$ and edges corresponding to the generating equations in $\RR$. In Remark \ref{no_normal_forms}, we defined a functor $(\widehat{-},\widetilde{-})$ from $\Gamma_n^{\GG,\RR}$ to the groupoid whose objects are pairs $(C,x)$ where $C$ is a computad and $x\in\Tree_n^{\GG,\RR}(C)$ is a term and whose morphisms $(C,x)\to (D,y)$ are isomorphisms $\varphi:C\to D$ such that $\varphi(x)\stackrel{\RR}{=}y$. This functor is easily seen to be full, so the group of automorphisms of $(\hat{x},\tilde{x})$ is a quotient of the group of automorphisms of $x$. In order to show that the category of computads for this $n$-globular operad is a presheaf category, it is enough to show that the former is trivial. In Remark \ref{no_normal_forms} we showed how one can use normal forms to do this.

On the other hand, one can consider the relation on the arrows of $\Gamma_n^{\GG,\RR}$ obtained by declaring that the application of a generating equation at any location in a tree commutes with the application of another equation in a disjoint location. The functor $(\widehat{-},\widetilde{-})$ respects this relation. If the automorphism groups of the quotient of $\Gamma_n^{\GG,\RR}$ by this relation are trivial, then the associated category of computads is a preseheaf category. It is an interesting question whether this is true for $n$-sesquicategories.

More generally, if one can describe a relation on the arrows of $\Gamma_n^{\GG,\RR}$ which is preserved by $(\widehat{-},\widetilde{-})$ and such that the automorphism groups of the quotient of $\Gamma_n^{\GG,\RR}$ by this relation are trivial, this proves that the associated category of computads is a preseheaf category.
 
\end{remark}

\section{String diagrams for $n$-sesquicategories}

Let $C$ be an $n$-computad for $T_n^{\sd}$. In this section we explain how to associate a \textbf{$C$-labelled string diagram} to each morphism in $F_n(C)$. This is an extremely useful graphical notation for describing composites and performing computations in $n$-sesquicategories. We will in the future extend this to semistrict $n$-categories by adding interchangers, which will allow us to apply in this more general context the techniques used in \cite{araujo_thesis},\cite{fibrations}, \cite{adj3} and \cite{adj4} (used there in the context of strict $3$ and $4$-categories). 

The essential ingredient here is the theory of normal forms, which will allow us to describe the graphical notation for a term by induction on the $(\OO_n,C)$-labelled tree corresponding to its normal form.

It is enough to decribe the unlabelled diagrams corresponding to morphisms in $F_n(\one)$, since $C$-labelled diagrams are then obtained by simply adding labels at appropriate places. 

So let $w\in F_n(\one)_k$. We proceed by induction on $k$ and the height of the normal form of $w$. When $k=0$, the morphism $w$ simply consists of the unique $0$-cell in $\one$. The associated diagram is just a point. In general, when $k$ is odd (resp. even) we depict a generating $k$-cell $w\in \one_k$ by drawing the $(k-1)$-diagram corresponding to its source on the left (resp. top), the one corresponding to its target on the right (resp. bottom) and then forming  a double cone on this disjoint union. We denote this double cone by drawing the cone point in the middle and curves connecting each cell in the source and target diagrams to the cone point. We need to distinguish lines which correspond to cells of different codimension, which we can do by using different thickness, transparency, dashing or any other method.

Now suppose $w=(x\to\circ_{i,j}\leftarrow y)$ is in normal form. By induction, we already know how to draw the diagrams associated to $x$ and $y$ and the diagram $\circ_{i,j}$ determines how we should compose these two pictures to obtain the picture for $w$. Finally, if $w=(x\to u_k)$ in normal form, then we draw two copies of $x$ and we draw lines connecting generators, again using a different notation for generators of different codimension.

We now give some examples. First it is useful to recall fom \cite{sesquicat} the graphical notation for the generators $\circ_{i,j}$. We include here the pictures for $i,j\leq 4$.

\begin{center}\begin{tabular}{lllll}

$\circ_{1,1}=\includegraphics[align=c]{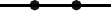}$ & $\circ_{1,2}=\includegraphics[align=c]{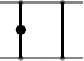}$ & $\circ_{2,1}=\includegraphics[align=c]{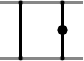}$ & $\circ_{2,2}=\includegraphics[align=c]{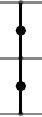}$ & $\circ_{1,3}=\includegraphics[align=c]{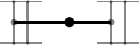}$
 
\end{tabular}\end{center}

\begin{center}\begin{tabular}{llll}

  $\circ_{3,1}=\includegraphics[align=c]{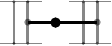}$ & $\circ_{2,3}=\includegraphics[align=c]{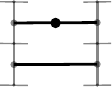}$ & $\circ_{3,2}=\includegraphics[align=c]{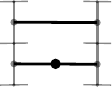}$ & $\circ_{3,3}=\includegraphics[align=c]{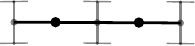}$
 
\end{tabular}\end{center}

\begin{center}\begin{tabular}{lllll}

  $\circ_{1,4}=\includegraphics[align=c]{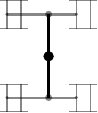}$ & $\circ_{4,1}=\includegraphics[align=c]{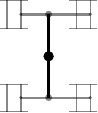}$ & $\circ_{2,4}=\includegraphics[align=c]{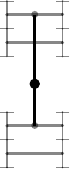}$ & $\circ_{4,2}=\includegraphics[align=c]{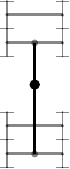}$ & $\circ_{3,4}=\includegraphics[align=c]{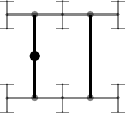}$
 
\end{tabular}\end{center}

\begin{center}\begin{tabular}{ll}

  $\circ_{4,3}=\includegraphics[align=c]{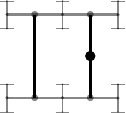}$ & $\circ_{4,4}=\includegraphics[align=c]{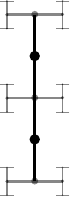}$
 
\end{tabular}\end{center}

Now we can move on to examples of unlabelled cells and diagrams. The unique $1$-cell in $\one$ is denoted $\includegraphics[align=c]{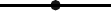}.$ The diagram $*\to u_1$ is denoted $\includegraphics[align=c]{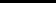}.$ The diagram $\includegraphics[align=c,scale=.5]{figures/1cell.pdf}\to\circ_{1,1}\leftarrow \includegraphics[align=c,scale=.5]{figures/1cell.pdf}$ is denoted $\includegraphics[align=c]{figures/circ11.pdf}.$ Here are some $2$-cells, with their source and target $1$-diagrams.

\begin{center}\begin{tabular}{llclc}
 
$\includegraphics[align=c]{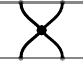}$ & $:$ & $\includegraphics[align=c]{figures/circ11.pdf}$ & $\to$ & $\includegraphics[align=c]{figures/circ11.pdf}$

\\

\\

$\includegraphics[align=c]{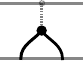}$ & $:$ & $\includegraphics[align=c]{figures/u1.pdf}$ & $\to$ & $\includegraphics[align=c]{figures/circ11.pdf}$

\\

\\

$\includegraphics[align=c]{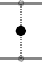}$ & $:$ & $\includegraphics[align=c]{figures/u1.pdf}$ & $\to$ & $\includegraphics[align=c]{figures/u1.pdf}$

\end{tabular}\end{center}

The $2$-diagrams $*\to u_1\to u_2$ and $\includegraphics[align=c,scale=.5]{figures/1cell.pdf}\to u_2$ are denoted \begin{center}\begin{tabular}{lll}

$\includegraphics[align=c]{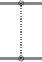}:\includegraphics[align=c]{figures/u1.pdf}\to\includegraphics[align=c]{figures/u1.pdf}$ & and & $\includegraphics[align=c]{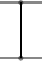}: \includegraphics[align=c]{figures/1cell.pdf}\to\includegraphics[align=c]{figures/1cell.pdf}.$
                                                                                                                                 
                                                                                                                    \end{tabular}
                                                                                                                                                                                                            \end{center}

Here is a $2$-diagram in normal form, then its normal form where we replace each generator by its picture and finally the picture of the diagram itself.

\[\vcenter{\vbox{\xymatrixcolsep{.6pc}\xymatrixrowsep{1pc}\xymatrix{ & & \circ_{2,2} & & \\ & \circ_{2,1}\ar[ru] & & \circ_{1,2}\ar[lu] & \\ \includegraphics[align=c,scale=.5]{figures/2cell_1.pdf}\ar[ru] & & \includegraphics[align=c,scale=.5]{figures/1cell.pdf}\ar[lu]\text{   }\includegraphics[align=c,scale=.5]{figures/1cell.pdf}\ar[ru] & & \includegraphics[align=c,scale=.5]{figures/2cell_1.pdf}\ar[lu]}}}=\includegraphics[align=c]{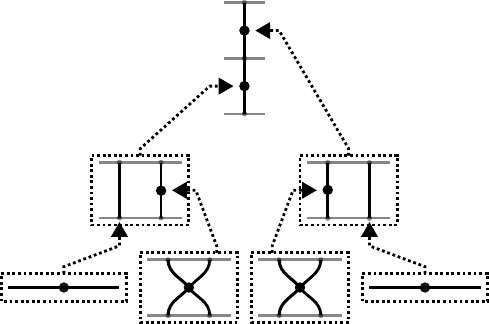}=\includegraphics[align=c]{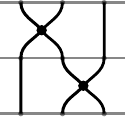}.\]

Here are some other 2-diagrams. $$\includegraphics[align=c]{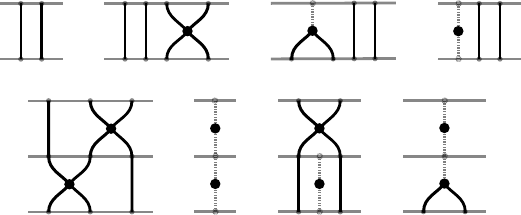}$$ Here are some $3$-cells with their source and target $2$-diagrams.

\begin{center}\begin{tabular}{llclc}
 
$\includegraphics[align=c]{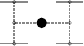}$ & $:$ & $\includegraphics[align=c]{figures/u2u1.pdf}$ & $\to$ & $\includegraphics[align=c]{figures/u2u1.pdf}$

\\

\\

$\includegraphics[align=c]{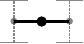}$ & $:$ & $\includegraphics[align=c]{figures/2cell_3.pdf}$ & $\to$ & $\includegraphics[align=c]{figures/2cell_3.pdf}$

\\

\\

$\includegraphics[align=c]{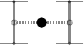}$ & $:$ & $\includegraphics[align=c]{figures/u2.pdf}$ & $\to$ & $\includegraphics[align=c]{figures/u2.pdf}$

\\

\\

$\includegraphics[align=c]{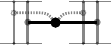}$ & $:$ & $\includegraphics[align=c]{figures/circ21.pdf}$ & $\to$ & $\includegraphics[align=c]{figures/circ21.pdf}$

\\

\\

$\includegraphics[align=c]{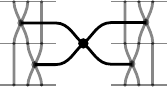}$ & $:$ & $\includegraphics[align=c]{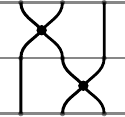}$ & $\to$ & $\includegraphics[align=c]{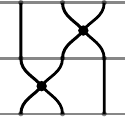}$

\end{tabular}\end{center}

Notice how the notation distinguishes the $3$-cell \includegraphics[align=c]{figures/3cell_4.pdf} above from the diagram $$\includegraphics[align=c]{figures/circ31.pdf}=(\includegraphics[align=c]{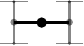}\to\circ_{3,1}\leftarrow\includegraphics[align=c]{figures/1cell.pdf}).$$ Here is a $4$-cell with its source and target $3$-diagrams.

$$\includegraphics[align=c]{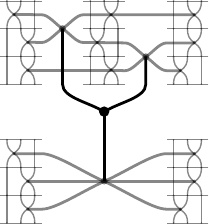}: \includegraphics[align=c]{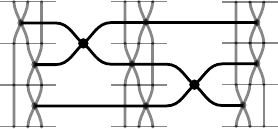}\to\includegraphics[align=c]{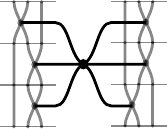}$$


\vspace{5mm}
\noindent
Manuel Ara\'{u}jo \\
manuel.araujo@tecnico.ulisboa.pt

\end{document}